\documentclass[12pt]{article}
\usepackage{amsmath,amssymb,amsthm}
\usepackage{graphicx,psfrag,epsf,float,color,subcaption}
\usepackage[caption = false]{subfig}
\usepackage{enumerate}
\usepackage{natbib}
\usepackage{url} 
\allowdisplaybreaks

\newcommand{\blind}{0}

\newtheorem{lem}{Lemma}
\newtheorem{thm}{Theorem}

\newtheorem{pro}{Proposition}
\theoremstyle{definition} \newtheorem{rem}{Remark}
\theoremstyle{definition} \newtheorem{ex}{Example}

\newcommand{\bfm}[1]{\ensuremath{\mathbf{#1}}}
\newcommand{\bdm}[1]{\ensuremath{\boldsymbol{#1}}}

\def \a {\bfm{a}}     \def \b {\bfm{b}}     
\def \d {\bfm{d}}          \def \f{\bfm{f}}
\def \g {\bfm{g}}

          \def \u{\bfm{u}}
                           \def \x{\bfm{x}}
\def \y {\bfm{y}}     

\def \A {\bfm{A}}     \def \B {\bfm{B}}     \def \C{\bfm{C}}
                           \def \F{\bfm{F}}
\def \G {\bfm{G}}     \def \H {\bfm{H}}     \def \I{\bfm{I}}
\def \J {\bfm{J}}          
\def \M {\bfm{M}}          
     \def \Q {\bfm{Q}}     \def \R{\bfm{R}}
\def \S {\bfm{S}}          
\def \V {\bfm{V}}          
\def \Y {\bfm{Y}}     

       \def \deltab   {\bdm{\delta}}
     \def \zetab    {\bdm{\zeta}}
         \def \thetab   {\bdm{\theta}}
        
      \def \mub      {\bdm{\mu}}

         \def \Sigmab   {\bdm{\Sigma}}

         \def \Omegab   {\bdm{\Omega}}

         \def \Sigmas   {{\Sigma}^{\ast}}

     \def\bB{\bfm B}

\def\bff{\bfm f}

\renewcommand{\hat}{\widehat}

\def \heps     {\hat{\heps}}

\newcommand{\norm}[1]{\lVert#1\rVert}
\newcommand{\bignorm}[1]{\Big\lVert#1\Big\rVert}
\newcommand{\fnorm}[1]{\lVert#1\rVert_F}
\newcommand{\bigfnorm}[1]{\Big\lVert#1\Big\rVert_F}
\newcommand{\snorm}[1]{\lVert#1\rVert_{\bdm{\Sigma}_{S}}}
\newcommand{\lonenorm}[1]{\lVert#1\rVert_1}

\newcommand{\maxnorm}[1]{\lVert#1\rVert_{\max}}

\newcommand{\hsigus}[1]{\hat{\Sigmab}_{u,S}^{#1}}

\DeclareMathOperator*{\argmin}{argmin}

\def \dd       {\mathrm{d}}
\def \E        {\mathrm{E}}
\def \cov      {\mathrm{cov}}
\def \var      {\mathrm{Var}}

\def \tr       {\mathrm{tr}}

\def \yt       {\y_{t}}
\def \ft       {\f_{t}}
\def \ut       {\u_{t}}

\def \Sc       {S^{c}}
\def \Sigmas   {\Sigmab_{S}}
\def \frateone {\frac{1}{\sqrt{T}}+\frac{1}{\sqrt{s}}}
\def \fratetwo {\frac{1}{\sqrt{T}}+\frac{1}{\sqrt{p}}}

\def \wt       {\tilde{\Sigmab}_u^{-1}}
\def \w        {{\Sigmab}_u^{-1}}
\def \GG       {\I_K+\H_1\B_S'\Sigmab_{u,S}^{-1}\B_S\H_1'}
\def \GGhat    {\I_K+\hat{\B}_1'(\hsigus{(1)})^{-1}\hat{\B}_1}

\begin{document}

\def\spacingset#1{\renewcommand{\baselinestretch}%
  {#1}\small\normalsize} \spacingset{1}


\if0\blind { \title{\bf Embracing the Blessing of Dimensionality in Factor
    Models} \author{ Quefeng Li, Guang Cheng, Jianqing Fan and Yuyan Wang}
  \maketitle
} \fi

\if1\blind
{
  \bigskip
  \bigskip
  \bigskip
  \begin{center}
    {\LARGE\bf Embracing the Blessing of Dimensionality in Factor Models}
  \end{center}
  \medskip
} \fi

\bigskip

\begin{abstract}
  Factor modeling is an essential tool for exploring intrinsic dependence structures among high-dimensional
  random variables. Much progress has been made for estimating the covariance matrix from a high-dimensional
  factor model. However, the blessing of dimensionality has not yet been fully embraced in the literature:
  much of the available data is often ignored in constructing covariance matrix estimates.  If our goal is to
  accurately estimate a covariance matrix of a set of targeted variables, shall we employ additional data,
  which are beyond the variables of interest, in the estimation? In this paper, we provide sufficient
  conditions for an affirmative answer, and further quantify its gain in terms of Fisher information and
  convergence rate. In fact, even an oracle-like result (as if all the factors were known) can be achieved
  when a sufficiently large number of variables is used. The idea of utilizing data as much as possible brings
  computational challenges. A divide-and-conquer algorithm is thus proposed to alleviate the computational
  burden, and also shown not to sacrifice any statistical accuracy in comparison with a pooled
  analysis. Simulation studies further confirm our advocacy for the use of full data, and demonstrate the
  effectiveness of the above algorithm. Our proposal is applied to a microarray data example that shows
  empirical benefits of using more data.
\end{abstract}

\noindent%
\vfill

\newpage
\spacingset{1.45} 

\section{Introduction}
\label{sec:1}
With the advance of modern information technology, it is now possible to track millions of variables or subjects simultaneously. To discover the relationship among them, the estimation of a high-dimensional
covariance matrix $\Sigmab$ has recently received a great deal of attention in the literature. Researchers proposed
various regularization methods to obtain consistent estimators of $\Sigmab$
\citep{BL08,LF09,RBLZ08,CL11,CTT10}. A key assumption for these regularization methods is that $\Sigmab$ is
sparse, i.e. many elements of $\Sigmab$ are small or exactly zero.

Different from such a sparsity condition,
factor analysis assumes that the intrinsic dependence is mainly driven by some common latent factors \citep{JW92}. For
example, in modeling stock returns, \cite{FF93} proposed the well-known Fama-French three-factor model. In the
factor model, $\Sigmab$ has spiked eigenvalues and dense entries. In the high dimensional
setting, there are many recent studies on the estimation of the covariance matrix based on the factor
model \citep{FFL08,FLM11,FLM13,BL12,BW13}, where the number of variables can be much larger than the number of observations.  

The interest of this paper is on the estimation of the covariance matrix for a certain set of variables using auxiliary data information.  In the literature, we use only the data information on the variables of interest.  In the data-rich environment today, substantially more amount of data information is indeed available, but is often ignored in statistical analysis. For example, we might be interested in understanding the covariance matrix of 50 stocks in a portfolio, yet the available data information are a time series of thousands of stocks.  Similarly,
an oncologist may wish to study the dependence or network structures among 100 genes that are significantly associated with
a certain cancer, yet she has expression data for over 20,000 genes from
the whole genome.  Can we benefit from using much more rich auxiliary data?

The answer to the above question is affirmative when a factor model is imposed.  Since the whole system is
driven by a few common factors, these common factors can be inferred more accurately from a much larger set of
data information \citep{FLM13}, which is indeed a ``blessing of dimensionality''.  A major contribution of
this paper is to characterize how much the estimation of the covariance matrix of interest and also common
factors can be improved by auxiliary data information (and under what conditions).

Consider the following factor model for all $p$ observable data  $\yt=(y_{1t},\ldots,y_{pt})'\in \mathbb{R}^p$ at time $t$:
\begin{equation}
  \label{eq:1}
  \yt=\B\ft+\ut, \quad t=1,\ldots,T,
\end{equation} where
$\ft \in\mathbb{R}^K$ is a $K$-dimensional vector of common factors, $\B=(\b'_1,\ldots,\b'_p)'\in \mathbb{R}^{p\times K}$ is a factor loading matrix with $\b_i\in \mathbb{R}^K$ being the factor loading of the $i$th variable on the latent factor $\bff_t$, and $\u_t$ is an idiosyncratic error vector. In the above model, $\y_t$ is the only observable variable, while $\bB$ is a matrix of unknown parameters, and $(\f_t, \u_t)$ are
latent  random variables.  Without loss of generality, we assume $\E(\f_t)=\E(\u_t)=\bdm{0}$ and $\f_t$ and
$\u_t$ are uncorrelated.  Then, the model implied covariance structure is
$$
\Sigmab =  \bB \cov(\bff_t) \bB' + \Sigmab_u,    
$$
where $\Sigmab=\E(\y_t\y_t')$ and $\Sigmab_u=\E(\u_t\u_t')$.
Observe that $\B$ and $\f_t$ are not individually
identifiable, since $\B\f_t=\B\H\H'\f_t$ for any orthogonal matrix $\H$. To this end, an identifiability
condition is imposed:
\begin{equation}
  \label{eq:2}
  \cov(\f_t)=\I_K \text{ and } \B'\Sigmab_u^{-1}\B \text{ is diagonal},
\end{equation}
which is a common assumption in the literature
\citep{BL12,BW13}.

Assume that we are only interested in a subset $S$ among a total of $p$ variables in model (\ref{eq:1}). We aim to obtain an efficient estimator of
\[
  \Sigmab_S=\B_S\B_S'+\Sigmab_{u,S},
\]
the covariance matrix of the $s$ variables in $S$, where $\B_S$ is the submatrix of $\B$ with row indices in
$S$ and $\Sigmab_{u,S}$ is the submatrix of $\Sigmab_u$ with row and column indices in $S$. As mentioned
above, the existing literature uses the following conventional method:
\begin{itemize}
\item Method 1: Use solely the $s$ variables in the set $S$ to estimate common factors $\f_t$, the loading matrix $\B_S$, the
  idiosyncratic matrix $\Sigmab_{u,S}$, and the covariance matrix $\Sigmab_S$.
\end{itemize}
This idea is apparently strongly influenced by the nonparametric estimation of the covariance matrix and  ignores a large portion of the available data in the other $p-s$ variables. An intuitively more efficient method is
\begin{itemize}
\item Method 2: Use all the $p$ variables to obtain estimators of $\f_t$, the loading matrix $\B$, the idiosyncratic matrix $\Sigmab_{u}$, and the entire covariance matrix $\Sigmab$, and then restrict them to the variables  of interest. This is the same as estimating $\f_t$ using all variables, and then estimating $\B_S$ and $\Sigmab_{u,S}$ based on the model \eqref{eq:1} and the subset $S$ with $\f_t$ being estimated (observed), and obtaining a plug-in estimator of $\Sigmab_S$.
\end{itemize}

We will show that Method 2 is more efficient than Method 1 in the estimation of $\f_t$ and $\Sigmab_S$ as more
auxiliary data information is incorporated. By treating common factor as an unknown parameter, we calculate
its Fisher information that grows with more data being utilized in Method 2. In this case, a more efficient
factor estimate can be obtained, e.g., through weighted principal component (WPC) method \citep{BW13}. The
advantage of factor estimation is further carried over to the estimation of $\Sigmab_S$ by Method 2 in terms
of its convergence rate. Moreover, if the number of total variables is sufficiently large, Method 2 is proven
to perform as well as an ``oracle method'', which observes all latent factors. This lends further support to
our aforementioned claim of ``blessing of dimensionality.'' Such a best possible rate improvement is new to
the existing literature, and counted as another contribution of this paper. All these conclusions hold when
the number of factors $K$ is assumed to be fixed and known, while $s$, $p$ and $T$ all tend to infinity.

The idea of utilizing data as much as possible brings computational challenges. Fortunately, we observe that all the $p$ variables are controlled by the same group of latent factors. Having said that, we can actually \emph{split $p$ variables} into smaller groups, and then utilize each group to estimate latent factors. The final factor estimate is obtained by averaging over these repeatedly estimated factors. Obviously, this divide-and-conquer algorithm can be implemented in a parallel computing environment, and thus produces factor estimators in a much more efficient way. On the other hand, our theory illustrates that this new method performs as well as the ``pooled analysis'', where we run the method over the whole dataset. Simulation studies further demonstrate the boosted computational speed and satisfactory statistical performance.

The rest of the paper is organized as follows. We compare the Fisher information of the factors by the two
methods in Section \ref{sec:2}. 
Section \ref{sec:3} describes the WPC method. As a main result, the convergence
rates of different estimators of $\Sigmab_S$ are further compared in Section
\ref{sec:4} under various norms. Section \ref{sec:5} introduces the
divide-and-conquer method for accelerating computation, while Section \ref{sec:6}
presents all simulation results. Section 7 gives a microarray data example to
illustrate our proposal. All technical proofs are delegated to the Appendix.

For any vector $\a$, let $\a_S$ denote a sub-vector of $\a$ with indices in $S$. Denote $\norm{\a}$ the
Euclidean norm of $\a$. For a symmetric matrix $\A\in \mathbb{R}^{d\times d}$, let $\A_{I,J}$ be the submatrix
of $\A$ with row and column indices in $I$ and $J$, respectively. We write $\A_S$ for $\A_{S,S}$ for
simplicity. Let $\lambda_j(\A)$ be the $j$th largest eigenvalue of $\A$. Denote
$\norm{\A}=\max \{|\lambda_1(\A)|, |\lambda_d(\A)|\}$ the operator norm of $\A$,
$\maxnorm{\A}=\max_{ij} |a_{ij}|$ the max-norm of $\A$, where $a_{ij}$ is the $(i,j)$-th entry of $\A$,
$\lonenorm{\A}=\max_i \sum_{j=1}^d |a_{ij}|$ the $L_1$ norm of $\A$, $\fnorm{\A}=\sqrt{\tr(\A'\A)}$ the
Frobenius norm of $\A$, and $\norm{\A}_{\M}=d^{-1/2} \fnorm{\M^{-1/2}\A\M^{-1/2}}$ the relative norm of $\A$
to $\M$, where the weight matrix $\M$ is assumed to be positive definite. For a non-square matrix $\C$, let
$\C_S$ be the submatrix of $\C$ with row indices in $S$.

\section{Fisher Information of Common Factor}
\label{sec:2}

In this section, we treat the vector of common factors as a fixed unknown parameter, and compute its Fisher
information matrices based on Method 1 and Method 2. In the computation, the loading matrix $\B$ is treated as
deterministic in Proposition 2. In Proposition 3, the Fisher information is computed for each given $\B$ and
then averaged over $\bB$ by regarding it as a realization of a chance process, which bypasses the block
diagonal assumption needed without taking average over $\B$. In other sections, we adopt the convention
regarding the factors as random and $\B$ as fixed. We start by calculating the Fisher information of
$\thetab_t:=\B\f_t$, which serves as an intermediate step in obtaining that for $\f_t$. For simplicity of
notation, time $t$ is suppressed in $(\y_t, \f_t, \u_t, \thetab_t)$ so that it becomes $(\y,\f, \u,\thetab)$
in this section.

Given a general density function of $\y$, denoted as $h(\y; \thetab)$, the Fisher information of $\thetab$
contained in full data is given by
\begin{equation*}
  I_p(\thetab)=\E \left[\left(\frac{\partial \log h(\y;\thetab)}{\partial \thetab} \right)\left(\frac{\partial \log
        h(\y;\thetab)}{\partial \thetab} \right)' \right].
\end{equation*}
When only data in $S$ is used, the Fisher information of $\thetab_S$ is given by
\begin{equation*}
  I_S(\thetab_S)=\E \left[\left(\frac{\partial \log h_S(\y_S;\thetab_S)}{\partial \thetab_{S}} \right)\left(\frac{\partial \log h_S(\y_S;\thetab_S)}{\partial \thetab_S} \right)' \right],
\end{equation*}
where $h_S$ is the marginal density of $\y_S$ for the target set of variable $S$. Our first proposition shows that $\{I_p(\thetab) \}_S$, the submatrix of $I_p(\thetab)$ restricted on $S$,
dominates $I_S (\thetab_S)$ under a mild condition.

\begin{pro}
  \label{pro:1}
  If $h(\y;\thetab)=h(\y-\thetab)$ and the density function $h(\y-\thetab)$ satisfies the following
  regularity condition: 
  \begin{equation}
    \label{eq:3}
    \nabla_{\y_S} \int h(\y_S-\thetab_S,\y_{\Sc}-\thetab_{\Sc}) \dd\y_{\Sc}= \int \nabla_{\y_S}
    h(\y_S-\thetab_S,\y_{\Sc}-\thetab_{\Sc}) \dd\y_{\Sc},
  \end{equation}
  then $\{I_p(\thetab)\}_S\succeq I_S(\thetab_S)$ in the sense that
  $\{I_p(\thetab)\}_S-I_S(\thetab_S)$ is positive semi-definite.
\end{pro}

The regularity condition \eqref{eq:3} is fairly mild, as illustrated in the following examples.
\begin{ex}
  \label{ex:1}
  In model \ref{eq:1}, if $\u_S$ and $\u_{\Sc}$ are independent, then \eqref{eq:3} holds.
\end{ex}
\begin{ex}
  \label{ex:2}
  If $\y$ follows an elliptical distribution that
  \begin{equation*}
    h(\y;\thetab)\propto g((\y-\thetab)'\Sigmab^{-1}(\y-\thetab)),
  \end{equation*}
  where the mapping function $g(t): [0,\infty) \to [0,\infty)$ satisfies that $|g'(t)|\leq c g(t)$ for some
  positive constant $c$, and $\E|\y|<\infty$, then \eqref{eq:3} holds. Example \ref{ex:2} includes some
  commonly used multivariate distributions as its special cases, e.g. the multivariate normal distribution and
  the multivariate $t$-distribution with degrees of freedom greater than 1. The proof is given in the Appendix
  Section \ref{sec:A2}.
\end{ex}

We next compute the Fisher information of $\f$ based on the full data set, denoted as $I(\f)$, and the partial
data set restricted on S, denoted as $I_S(\f)$. This can be done easily by noting that
$I(\f)=\B'I_p(\thetab)\B$. Indeed, the WPC estimators used in Methods 1 and 2 achieve such efficiency since
their asymptotic variances are proven to be the inverse of $I(\f)$ and $I_S(\f)$, respectively; see Remark
\ref{rem:1}.

Proposition \ref{pro:2} shows that $I(\f)$ dominates $I_S(\f)$, if $I_p(\thetab)$ is block-diagonal, i.e., $\{I_p(\thetab)\}_{S,S^c}=\bdm{0}$. Hence, common factors can be estimated more efficiently using additional data $\y_{S^c}$. The above block-diagonal condition implies that the idiosyncratic error of additional variables cannot be confounded with that of the variables-of-interest. For example, if $\u$ is normal, then $\{I_p(\thetab)\}_{S,S^c}=\bdm{0}$ indeed requires that $\u_S$ is independent of $\u_{S^c}$.

\begin{pro}
  \label{pro:2}
  Under condition (\ref{eq:3}), if $\{I_p(\thetab)\}_{S,S^c}=\bdm{0}$, $I(\f)\succeq I_S(\f)$.
\end{pro}

So far we treat $\B$ as being deterministic. Rather, Proposition \ref{pro:3} regards $\{\b_i\}$ as a
realization of a chance process. Under this assumption, the expectation of $I(\f)$ over $\B$ is shown to
always dominate that of $I_S(\f)$. In other words, we can claim that averaging over loading matrices, a larger
dataset contains more information about the unknown factors.

\begin{pro}
  \label{pro:3}
  If $\{\b_i \}_{i=1}^p$ are i.i.d. random loadings with $\E(\b_i)=\bdm{0}$ and (\ref{eq:3}) holds, then
  $\E[I(\f)]\succeq \E[I_S(\f)]$, where the expectation is taken with respect to the
  distribution of $\B$.
\end{pro}


\section{Efficient Estimation of Common Factor}
\label{sec:3}
In this section, we construct an efficient estimator of the common factors by showing that its asymptotic variance is exactly the inverse of its Fisher information. This together with the arguments in Section \ref{sec:2} enables us to draw a conclusion that using more data results in a more efficient factor estimator with a smaller asymptotic variance.

From a least-squares perspective, when the loading matrix $\B$ is known, $\f_t$
can be estimated by the weighted least-squares:
$\argmin_{\f_t\in \mathbb{R}^K}
\sum_{t=1}^T(\y_t-\B\f_t)'\Sigmab_u^{-1}(\y_t-\B\f_t)$.
In the high-dimensional setting ($p \gg T$), we assume $\Sigmab_u$ is a sparse
matrix and define its sparsity measurement as
\begin{equation}
  \label{eq:4}
  m_p=\max_{i\leq p} \sum_{j\neq i} I(\sigma_{u,ij}\neq 0),  \text{ where } \sigma_{u,ij} \text{ is the } (i,j)\text{-th
    entry of } \Sigmab_u.
\end{equation}
In particular, we assume the following sparsity condition
\begin{equation}
  \label{eq:5}
  m_p= o \left(\min \left\{\frac{1}{p^{1/4}} \sqrt{\frac{T}{\log p}}, p^{1/4} \right\} \right) \text{ and }
  \sum_{i=1}^p \sum_{j\neq i} I(\sigma_{u,ij} \neq 0) =O(p).
\end{equation}
Now, we propose to solve the following constrained weighted least-squares problem:
\begin{equation}
  \label{eq:6}
  \begin{split}
    (\hat{\B},\hat{\f}_1,\ldots, \hat{\f}_T) &=\argmin_{\B,\f_t}~ \sum_{t=1}^T
    (\y_t-\B\f_t)'\tilde{\Sigmab}_u^{-1}(\y_t-\B\f_t),\\
    \text{subject to } &\frac{1}{T} \sum_{t=1}^T {\f}_t{\f}_t'=\I_K; ~
    {\B}'\tilde{\Sigmab}_u^{-1}{\B} \text{
      is diagonal,}
  \end{split}
\end{equation}
where $\tilde{\Sigmab}_u$ is a regularized estimator of $\Sigmab_u$ to be discussed later. The above constraint is a sample analog of the identifiability condition (\ref{eq:2}). The involvement of the weight $\tilde{\Sigmab}_u^{-1}$
is to account for the heterogeneity among the data and leads to more efficient estimation of $(\B,\ft)$
\citep{CH12,BW13}.

Indeed, an initial estimator $\tilde{\Sigmab}_u$ of the idiosyncratic matrix $\Sigmab_u$ is needed for solving the
constrained weighted least-squares problem. We propose to obtain such an estimator by the following procedure, which is in
the same spirit as the estimation of the idiosyncratic matrix in the POET method \citep{FLM13}. Let
$\S_y= {T}^{-1} \sum_{t=1}^T (\y_t-\bar{\y})(\y_t-\bar{\y})'$ be the sample covariance of $\y$ and
$\{(\lambda_i,\zetab_i) \}_{i=1}^p$ be eigen-pairs of $\S_y$ with $\lambda_1\geq \lambda_2\geq
\cdots\geq\lambda_p$.
Denote $\R=\S_y-\sum_{i=1}^K \lambda_i\zetab_i\zetab_i'$. We estimate $\Sigmab_u$ by $\tilde{\Sigmab}_u$, whose $(i,j)$-th
entry
\begin{equation*}
  \hat{\sigma}_{u,ij}=
  \begin{cases}
    r_{ii}, & \text{ for } i=j,\\
    s_{ij}(r_{ij}), & \text{ for } i\neq j,
  \end{cases}
  \quad \text{ where } \R=(r_{ij}),
\end{equation*}
$s_{ij}(r_{ij})$ is a general entry-wise thresholding function \citep{AF01} such that
$s_{ij}(z)=0$ if $|z|\leq\tau_{ij}$ and $|s_{ij}(z)-z|\leq \tau_{ij}$ for
$|z|>\tau_{ij}$. In our paper, we choose hard-thresholding even though SCAD
\citep{FL01} and MCP \citep{ZHA10} are also applicable. We specify the entry-wise
thresholding level as
\begin{equation}
  \label{eq:7}
  \tau_{ij}(p)=C \sqrt{r_{ii}r_{jj}} \omega(p), \text{ where } \omega(p)= \sqrt{\frac{\log p}{T}}+\frac{1}{\sqrt{p}},
\end{equation}
and $C$ is a constant chosen by cross-validation. The thresholding parameter $C\omega(p)$ is applied to the correlation matrix. This is similar to the adaptive thresholding estimator for a general
covariance matrix \citep{RLZ09}, where the entry-wise thresholding level depends on $p$.

With $\tilde{\Sigmab}_u$ being the thresholding estimator described above, the constrained weighted least-squares problem \eqref{eq:6} can be solved by the weighted principal component (WPC) method. The solution is
given by
\begin{equation}
  \label{eq:8}
  \hat{\F}=(\hat{\f}_1,\ldots,\hat{\f}_T)' \text{ and } \hat{\B}=T^{-1}\Y\hat{\F},
\end{equation}
where $\Y=(\y_1,\ldots,\y_T)$ and the columns of $\hat{\F}$ are the eigenvectors
corresponding to the largest $K$ eigenvalues of the $T\times T$ matrix
$\sqrt{T}\Y'\tilde{\Sigmab}_u^{-1}\Y$ \citep{BW13}.

In the following, we give a result showing that the WPC estimator is asymptotically efficient. Indeed,
\cite{BW13} derive the asymptotic normality of $\hat{\f}_t$ under the following conditions:
\begin{enumerate}
\item[(i)] All eigenvalues of $\B'\B/p$ are bounded away from
  zero and infinity as $p\to \infty$;

\item[(ii)] There exists a $K\times K$ diagonal matrix $\Q$ such that
  ${\B'\Sigmab_u^{-1}\B}/{p}\to \Q$. In addition, the diagonal elements of $\Q$
  are distinct and bounded away from infinity.

\item[(iii)] For each fixed $t\leq T$,
  $(\B'\Sigmab_u^{-1}\B)^{-1/2} \B'\Sigmab_u^{-1}\u_t \xrightarrow{d} N(\bdm{0},\I_K)$, as $p\to \infty$,
\end{enumerate}
together with the sparsity assumption \eqref{eq:5}, and some additional regularity conditions given in Section
\ref{sec:A1}. When $\sqrt{p}\log p=o(T)$, it is shown that
\begin{equation}
  \label{eq:9}
  \sqrt{p}(\hat{\f}_t-\H\f_t)\xrightarrow{D} N(\bdm{0}, \Q^{-1}),
\end{equation}
where $\H$ is a specific rotation matrix given by
\begin{equation}
  \label{eq:10}
  \H=\hat{\V}^{-1}\hat{\F}'\F\B'\tilde{\Sigmab}_u^{-1}\B/T,
\end{equation}
and $\hat{\V}$ is a $K\times K$ diagonal matrix of the largest $K$ eigenvalues of
$\Y'\tilde{\Sigmab}_u^{-1}\Y/T$. The rotation matrix $\H$ is introduced here so that $\H\f_t$
  is an identifiable quantity from the data. See more discussion about the identifiability in Remark
  \ref{rem:2}.

Condition (i) is a ``pervasive condition'' requiring that the common factors affect a non-negligible fraction
of subjects. This is a common assumption for the principal components based methods \citep{FLM11,BW13}. In
condition (ii), $\B'\Sigmab_u^{-1}\B$ is indeed the Fisher information (under Gaussian errors) contained in $p$ variables, while the limit $\Q$ can be viewed as an average information for each variable. Hence, the asymptotic normality in (\ref{eq:9})   shows that $\hat{\f}_t$ is efficient as its asymptotic variance attains the inverse of the (averaged) Fisher information.

\begin{rem}
  \label{rem:1}
  The results in Section \ref{sec:2} together with \eqref{eq:9} imply that Method 2 is in general better than
  Method 1 in the estimation of common factors. To explain why, we consider two different cases here. When $p$
  is an order of magnitude larger than $s$, where $s$ is the number of variables of interest. Method 2
  produces a better estimator of factors with a faster convergence rate. Even when $p$ and $s$ diverge at the
  same speed, the factor estimator based on Method 2 is shown to possess a smaller asymptotic variance, as
  long as $\Sigmab_{u,S,S^c}=\bdm{0}$. Recall that $\B'\Sigmab_u^{-1}\B=I(\f)$ and
  $\B_S'\Sigmab_{u,S}^{-1}\B_S=I_S(\f)$ under Gaussian errors, and they also correspond to the inverse of the
  asymptotic variance given by Methods 1 and 2, respectively. Then, Proposition 2 implies that Method 2 has a
  smaller asymptotic variance, if $\Sigmab_{u,S,S^c}=\bdm{0}$. Alternatively, if $\B$ is treated as being
  random, Proposition \ref{pro:3} immediately implies that
  $E(\B_S'\Sigmab_{u,S}^{-1}\B_S)\succeq E(\B'\Sigmab^{-1}\B)$.  Therefore, even without the block diagonal
  assumption, Method 2 produces a more efficient factor estimate on average.

\end{rem}

\section{Covariance Matrix Estimation}
\label{sec:4}
One primary goal in this paper is to obtain an accurate estimator of the covariance matrix $\Sigmas=\E(\y_S\y'_S)$ for the variables-of-interest. In this section, we compare three different estimation methods, namely Methods 1, 2 and Oracle Method, in terms of their rates of convergence (under various norms). Obviously, these rates depend on how accurately the realized factors are estimated as demonstrated later.

Below we describe these three methods in full details.
\begin{itemize}
\item Method 1:
  \begin{itemize}
  \item [i.] Use solely the data in the subset $S$ to obtain estimators of
    the realized factors $\hat{\F}^{(1)}$ and the loading matrix $\hat{\B}_1=T^{-1}\Y_S\hat{\F}^{(1)}$ based on (\ref{eq:8});

  \item [ii.] Let $(\hat{\f}_t^{(1)})'$ be the $t$-th row of $\hat{\F}^{(1)}$,
    $(\hat{\b}_i^{(1)})'$ be the $i$-th row of $\hat{\B}_1$,
    $\hat{u}_{it}=y_{it}-(\hat{\b}_i^{(1)})'\hat{\f}^{(1)}_t$, and
    $\hat{\sigma}_{ij}=\frac{1}{T} \sum_{t=1}^T \hat{u}_{it}\hat{u}_{jt}$.  The
    $(i,j)$-th entry of the idiosyncratic matrix estimator
    $\hat{\Sigmab}^{(1)}_{u,S}$ of $\Sigmab_{u,S}$ is given by thresholding
    $\hat{\sigma}_{ij}$ at the level of $C\hat{\theta}_{ij}^{1/2} \omega(s)$,
    where $\omega(s)$ is defined in (\ref{eq:7}) and
    $\hat{\theta}_{ij}=\frac{1}{T} \sum_{t=1}^T
    (\hat{u}_{it}\hat{u}_{jt}-\hat{\sigma}_{ij})^2$;
  \item [iii.] The final estimator is given by
    $\hat{\Sigmab}_S^{(1)}=\hat{\B}_1\hat{\B}_1'+ \hat{\Sigmab}_{u,S}^{(1)}$.

  \end{itemize}
\item Method 2:
  \begin{itemize}
  \item [i.] Use all $p$ variables to obtain the estimate
    $\hat{\F}^{(2)}$ as given in (\ref{eq:8}) for the realized factors
    and then estimate the loading
    $\B_S$ by $\hat{\B}_2=T^{-1}\Y_S\hat{\F}^{(2)}$;
  \item [ii.] Follow the same procedure as in Method 1 to obtain the estimator
    $\hat{\Sigmab}_{u,S}^{(2)}$ but based on $\hat{\F}^{(2)}$ and $\hat{\B}_2$;
  \item [iii.] The final estimator is given by
    $\hat{\Sigmab}_S^{(2)}=\hat{\B}_2\hat{\B}_2'+ \hat{\Sigmab}_{u,S}^{(2)}$.
  \end{itemize}
\item Oracle Method:
  \begin{itemize}
  \item [i.] Estimate the loading by $\hat{\B}_o=T^{-1} \Y_S\F$, where
    $\F=(\f_1,\ldots, \f_T)'$ are the true factors.
  \item [ii.] The idiosyncratic matrix estimator $\hat{\Sigmab}_{u,S}^o$ is given
    by the same procedure as in Method 1, with $\hat{\b}_i^{(1)}$ and
    $\hat{\f}_t^{(1)}$ being replaced by $\hat{\b}_i^o$ and $\f_t$, respectively.
  \item[iii.] The final estimator is given by
    $\hat{\Sigmab}_S^{o}=\hat{\B}_o\hat{\B}_o'+ \hat{\Sigmab}_{u,S}^{o}$.
  \end{itemize}
\end{itemize}

Theorem \ref{thm:1} depicts the estimation accuracy of $\Sigmab_S$ by the above three methods with respect to the following measurements:
\begin{align*}
  \snorm{\hat{\Sigmab}_S-\Sigmab_S}, \quad
  \maxnorm{\hat{\Sigmab}_S-\Sigmab_S}, \quad \norm{\hat{\Sigmab}_S^{-1}-\Sigmab_S^{-1}},
\end{align*}
where
$\snorm{\hat{\Sigmab}_S-\Sigmab_S} = p^{-1/2} \| \Sigmab_S^{-1/2} \hat{\Sigmab}_S \Sigmab_S^{-1/2} - \I_S\|_F$
is a norm of the relative errors. Note that the results of \cite{FLM13} can not be directly used here since we
employ the {\em weighted} principal component analysis to estimate the unobserved factors.  This is expected
to be more accurate than the ordinary principal component analysis, as shown in \cite{BW13}.  Indeed, the
technical proofs for our results are technically more involved than those in \cite{FLM13}. 

We assume that $s$ is much less than $p$, i.e., $s=o(p)$, but both tend to
infinity. Under the pervasive condition (i), $\norm{\Sigmab_S}\geq cs$ and
therefore diverges. For this reason, we consider the relative norm
$\snorm{\hat{\Sigmab}_S-\Sigmab_S}$, instead of
$\norm{\hat{\Sigmab}_S-\Sigmab_S}$, and the operator norm
$\norm{\hat{\Sigmab}_S^{-1}-\Sigmab_S^{-1}}$ for estimating the inverse. In
addition, we consider another element-wise max norm
$\maxnorm{\hat{\Sigmab}_S-\Sigmab_S}$. We show that if $p$ is large with respect
to $s$ and $T$, Method 2 performs as well as the Oracle Method, both of which
outperform Method 1. As a consequence, even if we are only interested in the
covariance matrix of a small subset of variables, we should use all the data to
estimate the common factors, which ultimately improves the estimation of
$\Sigmas$. In particular, we are able to specify an explicit regime of $(s,p)$
under which the improvements are substantial. However, when $s \asymp p$,
i.e. they are in the same order, using more data does not show as dramatic
improvements for estimating $\Sigmab_S$. This is expected and will be clearly
seen in the simulation section.

Before stating Theorem \ref{thm:1}, we need a few preliminary results: Lemmas
\ref{lem:1} -- \ref{lem:3}. Specifically, Lemma \ref{lem:1} presents the
uniform convergence rates of the factor estimates by Methods 1 and 2. Based
on that, Lemmas \ref{lem:2} and \ref{lem:3} further derive the estimation accuracy of factor loadings and idiosyncratic matrix by the three methods,
respectively. These results together lead to the estimation error rates of
$\Sigmab_S$ in Theorem \ref{thm:1} w.r.t. three measures defined above. Additional Lemmas supporting the proof are given in Appendix.  Again, these kinds of results can not be obtained directly from \cite{FLM13} due to our use of WPC.

\begin{lem}
  \label{lem:1}
  Suppose that conditions (i), (ii), the sparsity condition (\ref{eq:5}), and additional regularity conditions
  (iv)-(vii) in Section \ref{sec:A1} hold for both $s$ and $p$. If $\sqrt{p}\log p=o(T)$ and $T=o(s^2)$, then
  we have
  \begin{equation*}
    \max_{t\leq T} \norm{\hat{\f}^{(1)}_t-\H_1 \f_t}=O_P\left(\frac{1}{\sqrt{T}}+\frac{T^{1/4}}{\sqrt{s}}\right)
    \quad \text{ and } \quad
    \max_{t\leq T} \norm{\hat{\f}^{(2)}_t-\H_2 \f_t} =O_P\left(\frac{1}{\sqrt{T}}+\frac{T^{1/4}}{\sqrt{p}}\right),
  \end{equation*}
  where
  $\H_1=\hat{\V}_1^{-1}\hat{\F}^{(1)'}\F\B_S'\tilde{\Sigmab}_{u,S}^{-1}\B_S/T$,
  $\H_2=\hat{\V}_2^{-1}\hat{\F}^{(2)'}\F\B'\tilde{\Sigmab}_u^{-1}\B/T$,
  $\hat{\V}_1$ is the diagonal matrix of the largest $K$ eigenvalues of
  $\Y_S'\tilde{\Sigmab}^{-1}_{u,S}\Y_S/T$ and $\hat{\V}_2$ is the diagonal matrix
  of the largest $K$ eigenvalues of $\Y'\tilde{\Sigmab}^{-1}_u\Y/T$. 
\end{lem}
\begin{rem}
  \label{rem:2}
  $\H_1$ and $\H_2$ correspond to the rotation matrix $\H$ defined in (10) using Methods 1 and 2,
  respectively. Recall that $\F=(\f_1,\ldots, \f_T)'$, then
  $\H\f_t=T^{-1}\hat{\V}^{-1}\hat{\F}(\B\f_1,\ldots,\B\f_T)'\tilde{\Sigmab}_u^{-1} \B\f_t$. Note that $\H\f_t$
  only depends on quantities $\V^{-1}\hat{\F}$, $\tilde{\Sigmab}_u^{-1}$ and the {\em identifiable} component
  $\{\B\f_t \}_{t=1}^T$. Therefore, there is no identifiability issue regarding $\H\f_t$. In other words, even
  though $\f_t$ itself may not be identifiable, an identifiable rotation of $\f_t$ can be consistently
  estimated by $\hat{\f}_t$.
\end{rem}

\noindent Lemma \ref{lem:1} implies that Method 2 produces a better factor estimate if $$0.5<\gamma_s<1.5\leq\gamma_p<2,$$ by
representing $s$ and $p$ as $s\asymp T^{\gamma_s}$ and $p \asymp T^{\gamma_{p}}$.

It is not surprising that the estimation accuracy of
loading matrix also varies among these three methods as shown in Lemma \ref{lem:2} below.
\begin{lem}
  \label{lem:2}
  Under conditions of Lemma \ref{lem:1},
  \begin{align*}
    \max_{i\leq s} \norm{\hat{\b}_i^{(1)}-\H_1\b_i}
    &=O_P\left(w_1 \right), \text{ where } w_1:=\frac{1}{\sqrt{s}}+\sqrt{\frac{\log s}{T}},\\
    \max_{i\leq s} \norm{\hat{\b}_i^{(2)}-\H_2\b_i}
    &=O_P\left(w_2 \right), \text{ where } w_2:=\frac{1}{\sqrt{p}}+\sqrt{\frac{\log s}{T}},\\
    \max_{i\leq s} \norm{\hat{\b}_i^{o}-\b_i}
    &=O_P\left(w_o \right), \text{ where }
      w_o:=\sqrt{\frac{\log s}{T}}.
  \end{align*}
\end{lem}

Similarly, Lemma \ref{lem:2} indicates that Method 2 performs as well as the Oracle Method, both of which are better than Method 1, i.e., $w_2=w_o < w_1$, if
\begin{equation*}
  0.5<\gamma_s<1\leq\gamma_p<2,
\end{equation*}
by representing $s$ and $p$ in the order of $T$ as above. We remark that the extra terms $1/\sqrt{s}$ and $1/\sqrt{p}$ in $w_1$ and $w_2$ (in comparison with the oracle rate $w_o$) are due to the factor estimation. Another preliminary result regarding the estimation of the identifiable component
$\b_i'\f_t$ is given in Lemma \ref{lem:A1}.

Similar insights can be delivered from Lemma \ref{lem:3} on the estimation of $\Sigmab_{u,S}$.

\begin{lem}
  \label{lem:3}
  Under conditions of Lemma \ref{lem:1}, it holds that
  \begin{align*}
    \norm{\hat{\Sigmab}_{u,S}^{(1)}-\Sigmab_{u,S}}
    &=O_P\left(m_s w_1 \right)=\norm{(\hat{\Sigmab}_{u,S}^{(1)})^{-1}-\Sigmab_{u,S}^{-1}},\\
    \norm{\hat{\Sigmab}_{u,S}^{(2)}-\Sigmab_{u,S}}
    &=O_P\left(m_s w_2 \right)=\norm{(\hat{\Sigmab}_{u,S}^{(2)})^{-1}-\Sigmab_{u,S}^{-1}}, \\
    \norm{\hat{\Sigmab}_{u,S}^{o}-\Sigmab_{u,S}}
    &=O_P\left(m_s w_o\right)=\norm{(\hat{\Sigmab}_{u,S}^{o})^{-1}-\Sigmab_{u,S}^{-1}},
  \end{align*}
  where $m_s$ is defined as in (\ref{eq:4}) with $p$ being replaced by $s$.
\end{lem}

Now, we are ready to state our main result on the estimation of $\Sigmab_S$ based on the above preliminary
results. From Theorem \ref{thm:1}, it is easily seen that the comparison of the estimation accuracy of
$\Sigmab_S$ among three methods is solely determined by the relative magnitude of $w_o$, $w_1$ and
$w_2$. Therefore, we should use additional variables to estimate the factors if $p$ is much larger than $s$ in
the sense that $T/\log s=O(p)$ and $s\log s=o(T)$ (implying $w_2=w_o < w_1$).

\begin{thm}
  \label{thm:1}
  Under conditions of Lemma \ref{lem:1} , it holds that

  (1) For the relative norm, $\snorm{\hat{\Sigmab}_S^{(1)}-\Sigmab_S}=O_P\left(\sqrt{s}w_1^2+m_sw_1\right)$,
  $\snorm{\hat{\Sigmab}_S^{(2)}-\Sigmab_S}=O_P\left(\sqrt{s}w_2^2+m_sw_2 \right)$, and
  $\snorm{\hat{\Sigmab}_S^{o}-\Sigmab_S}=O_P\left(\sqrt{s}w_o^2+m_sw_o \right)$.

  (2) For the max-norm, $\maxnorm{\hat{\Sigmab}_S^{(1)}-\Sigmab_S}=O_P\left(w_1\right)$,
  $\maxnorm{\hat{\Sigmab}_S^{(2)}-\Sigmab_S}=O_P\left(w_2\right)$, and
  $\maxnorm{\hat{\Sigmab}_S^{o}-\Sigmab_S}=O_P\left(w_o\right)$.

  (3) For the operator norm of the inverse matrix, $\norm{(\hat{\Sigmab}_{S}^{(1)})^{-1}-\Sigmab^{-1}_{S}}=O_P\left(m_sw_1\right)$,
  $\norm{(\hat{\Sigmab}_{S}^{(2)})^{-1}-\Sigmab^{-1}_{S}}=O_P\left(m_sw_2 \right)$ and
  $\norm{(\hat{\Sigmab}_{S}^o)^{-1}-\Sigmab^{-1}_{S}}=O_P\left(m_sw_o \right)$.
\end{thm}

\begin{rem}
  \label{rem:3}
  So far, we assumed that the number of factors $K$ is fixed and known. A data driven choice
  of $K$ has been extensively studied in the econometrics literature, e.g., by \cite{BN02}, \cite{K10}. To estimate $K$, we can adopt the method by \cite{BN02} and propose a consistent estimator of $K$ (by
  allowing $p, T\to\infty$) as follows
  \begin{equation*}
    \hat{K}=\argmin_{0\leq k\leq N} \log \left\{\frac{1}{pT}  \fnorm{\Y-T^{-1}\Y
        \hat{\F}_k\hat{\F}_k'}^2\right\}+kg(p,T), 
  \end{equation*}
  where $N$ is a predefined upper bound, $\hat{\F}_k$ is a $T\times k$ matrix whose columns are $\sqrt{T}$
  times the eigenvectors corresponding to the largest $k$ eigenvalues of $\Y'\Y$, and $g(p,T)$ is a penalty
  function. Two examples suggested by \cite{BN02} are
  \begin{equation*}
    g(T,p)=\frac{p+T}{pT} \log \left(\frac{pT}{p+T} \right) \quad \text{ or } \quad
    g(T,p)=\frac{p+T}{pT} \log \left(\min\{p,T \} \right). 
  \end{equation*}
  Under our assumptions (i)-(x), all conditions required by theorem 2 of \cite{BN02} hold. Hence, their
  theorem implies that $P(\hat{K}=K)\to 1$. Then, conditioning on the event that $\{\hat{K}=K \}$, our theorem
  1 still holds by replacing $K$ with $\hat{K}$.  Other effective methods for selecting the number of factors
  include the eigen ratio method in \cite{LY12} and \cite{AH13}.
\end{rem}
\begin{rem}
  \label{rem:4}
  When $K$ grows with $p$ and $T$, Fan et al. (2013) gives the explicit
  dependence of the convergence rates on $K$ for their proposed POET
  estimator. By adopting
  their technique, we can obtain the following results: \\
  (1)
  $\snorm{\hat{\Sigmab}_S^{(1)}-\Sigmab_S}=O_P\left(K\sqrt{s}w_1^2+K^3m_sw_1\right)$,
  $\snorm{\hat{\Sigmab}_S^{(2)}-\Sigmab_S}=O_P\left(K\sqrt{s}w_2^2+K^3m_sw_2\right)$,
  $\snorm{\hat{\Sigmab}_S^{o}-\Sigmab_S}=O_P\left(K\sqrt{s}w_o^2+K^3m_sw_o\right)$;
  \\
  (2) $\maxnorm{\hat{\Sigmab}_S^{(1)}-\Sigmab_S}=O_P\left(K^3w_1\right)$,
  $\maxnorm{\hat{\Sigmab}_S^{(2)}-\Sigmab_S}=O_P\left(K^3w_2\right)$,
  $\maxnorm{\hat{\Sigmab}_S^{o}-\Sigmab_S}=O_P\left(K^3w_o\right)$;\\
  (3)
  $\norm{(\hat{\Sigmab}_{S}^{(1)})^{-1}-\Sigmab^{-1}_{S}}=O_P\left(K^3m_sw_1\right)$,
  $\norm{(\hat{\Sigmab}_{S}^{(2)})^{-1}-\Sigmab^{-1}_{S}}=O_P\left(K^3m_sw_2\right)$,
  $\norm{(\hat{\Sigmab}_{S}^o)^{-1}-\Sigmab^{-1}_{S}}=O_P\left(K^3m_sw_o\right)$. \\
  Again, the rate difference among three types of estimators only depends on
  $w_o$, $w_1$ and $w_2$. Therefore, the same conclusion (when $p$ is much larger
  than $s$, using additional variables improves the estimation of $\Sigmab_S$)
  can still be made even if $K$ diverges. As long as $K$ diverges
    in the rate that $K=o(\min\{1/(\sqrt{s}w_1^2), 1/(m_sw_1)^{1/3}\})$,
    $K=o(1/w_1^{1/3})$ or $K=o(1/(m_sw_1)^{1/3})$, the same blessing of
    dimensionality phenomena persist in terms of estimation consistency in
    relative norm, max norm, or operator norm of the inverse, respectively.
\end{rem} 

\section{Divide-and-Conquer Computing Method}
\label{sec:5}

As discussed previously, we prefer utilizing auxiliary data information as much as possible even we are only interested in the covariance matrix of some particular set of variables. But this can bring up heavy computational burden. This concern motivates a simple divide-and-conquer scheme that {\em splits all $p$ variables in $\Y$}. Without loss of generality, assume that $p$ rows of matrix $\Y$ can be evenly divided into $M$ groups with $p/M$ variables in each group. The $s$ variables of interest can possibly be assigned to different
groups.

{\bf Divide-and-Conquer Computation Scheme}
\begin{enumerate}\vspace{-0.1in}
\item In the $m$th group, obtain the initial estimator $\tilde{\Sigmab}_{u,m}$ by using the adaptive
  thresholding method as described in Section \ref{sec:3} based on the data in the $m$th group only.
\item Denote $\Y_m$ as the data vector corresponding to the variables in the
  $m$th group and let $\hat{\F}_m=(\hat{\f}_{m,1},\ldots,\hat{\f}_{m,T})'$, where
  its columns are the eigenvectors corresponding to the largest $K$ eigenvalues
  of the $T\times T$ matrix
  $\sqrt{T}\Y_m'\tilde{\Sigmab}_{u,m}^{-1}\Y_m$. The computation in the above
  two steps can be done in a parallel manner.
\item Average $\{\hat{\f}_{m,t}\}_{m=1}^M$ to obtain a single estimator of
  $\f_t$ as
  $$\bar{\f}_t=\frac{1}{M} \sum_{m=1}^M \hat{\f}_{m,t}.$$
  The loading matrix estimate is given by $\bar{\B}_S=T^{-1}\Y_S\bar{\F}$, where
  $\bar{\F}=(\bar{\f}_1,\ldots,\bar{\f}_T)'$.

\item The idiosyncratic matrix is estimated as follows. Let $\bar{\f}'_t$ be the
  $t$-th row of $\bar{\F}$ and $\bar{\b}'_i$ be the $i$th row of $\bar{\B}_S$. Let
  $\hat{u}_{it}=y_{it}-\bar{\b}'_i\bar{\f}_t$,
  $\hat{\sigma}_{ij}=T^{-1} \sum_{t=1}^T \hat{u}_{it}\hat{u}_{jt}$, and
  $\hat{\theta}_{ij}=T^{-1} \sum_{t=1}^T
  (\hat{u}_{it}\hat{u}_{jt}-\hat{\sigma}_{ij})^2$.
  The $(i,j)$-th entry of $\bar{\Sigmab}_{u,S}$ is given by
  thresholding $\hat{\sigma}_{ij}$ at the level of
  $C\hat{\theta}_{ij}^{1/2} \omega(s)$, where $\omega(s)$ is defined as in
  (\ref{eq:7}) with $p$ replaced by $s$.
\item The final estimator of the covariance matrix is given by
  \begin{equation*}
    \bar{\Sigmab}_S=\bar{\B}_S\bar{\B}_S'+ \bar{\Sigmab}_{u,S}.
  \end{equation*}
\end{enumerate}

We show that, if $M$ is fixed,
\begin{align*}
  \snorm{\bar{\Sigmab}_S-\Sigmab_S}&=O_P\left(\sqrt{s}w_2^2+m_sw_2 \right), \\
  \maxnorm{\bar{\Sigmab}_S-\Sigmab_S}&=O_P\left(w_2\right), \\
  \norm{(\bar{\Sigmab}_S)^{-1}-\Sigmab_S^{-1}}&=O_P\left(m_sw_2 \right).
\end{align*}
These rates match the rates of $\hat{\Sigmab}_S^{(2)}$ attained by Method 2, where all $p$ variables are
pooled together for the analysis. The proof is given in Appendix \ref{sec:A3}.
The simulation results in Section \ref{sec:6} further demonstrate that without
sacrificing the estimation accuracy, the divide-and-conquer method runs much
faster than Method 2. Therefore, the divide-and-conquer method is practically
useful when dealing with massive dataset.

The main computational cost of our method comes from taking the inverse of $\tilde{\Sigmab}_{u}$. For our
Method 2, where all $p$ variables are pooled together for the analysis, the computational complexity of the
inversion is $O(p^3)$. On the other hand, for the divide-and-conquer method, the corresponding estimator
$\tilde{\Sigmab}_{u,m}$ in the $m$-th group only needs a computational cost of $O((p/M)^3)$ to be
inverted. Then, the total computation complexity is $O(p^3/M^2)$. Hence, the computational speed can be
boosted by $M^2$-fold. Such a computational acceleration can also be observed from simulation study results in
Figure 1(d). Other operations like the eigen-decomposition on the $T\times T$ matrix
$\sqrt{T} \Y'\tilde{\Sigmab}_u^{-1}\Y$ do not have dominating computational cost, as we assume that $p$ is
much larger than $T$. When M grows too fast, the divide-and-conquer method may lose estimation efficiency
compared with the pooled analysis (Method 2). However, considering its boost of computation, the
divide-and-conquer method is practically useful when dealing with massive dataset.

\section{Simulations}
\label{sec:6}
We use simulated examples to compare the statistical performances of Methods 1, 2
and the Oracle Method. We fix the number of factors $K=3$ and repeat 100
simulations for each combination of $(s,p,T)$. The loading $\b_i$, the factor
$\f_t$ and the idiosyncratic error $\u_t$ are generated as follows:
\begin{itemize}
\item $\{\b_i \}_{i=1}^p$ are i.i.d. from $N_K(\bdm{0}, 5\I_K)$.
\item $\{\f_t \}_{t=1}^T$ are i.i.d. from $N_K(\bdm{0}, \I_K)$.
\item $\{\u_t \}_{t=1}^T$ are i.i.d. from $N_p(\bdm{0}, 50\I_p)$.
\end{itemize}
The observations $\{\y_t \}_{t=1}^T$ are generated from (\ref{eq:1}) using $\b_i$, $\f_t$ and $\u_t$ from the above. Tables \ref{table:1}-\ref{table:4} report the estimation errors of the factors, the loading matrices and the covariance-of-interest $\Sigmab_S$ in terms of different measurements.

We see from Tables \ref{table:1}  and \ref{table:2}  that when $s=50$ and
$p=1000, 2000$, Method 1 performs much worse than Method 2, for both $T=200$
and $T=400$.   However, when $s$
increases to $800$ with $p$ being the same, Tables \ref{table:3} and
\ref{table:4} show that the improvement of Method 2 over Method 1 is less
profound. This is expected as the set of interest already contains sufficiently
rich information to produce an accurate estimator for realized factors. In
general, we note that Method 2 is the most advantageous in the settings where $s$
is much smaller than $p$. In addition, from Tables 1-4, we can tell that Method 2
comes closer to the Oracle method as $p$
grows. 
In practice, we also observe that the WPC factor estimator performs better than the unweighted PC estimator
when $\u_t$ is heteroscedastic. Due to the space limit, we choose not to present the simulation results in
this model.

\begin{figure}[hbtp]
  \centering
  \begin{subfigure}[b]{0.5\textwidth}
    \centering
    \includegraphics[width=1.0\textwidth]{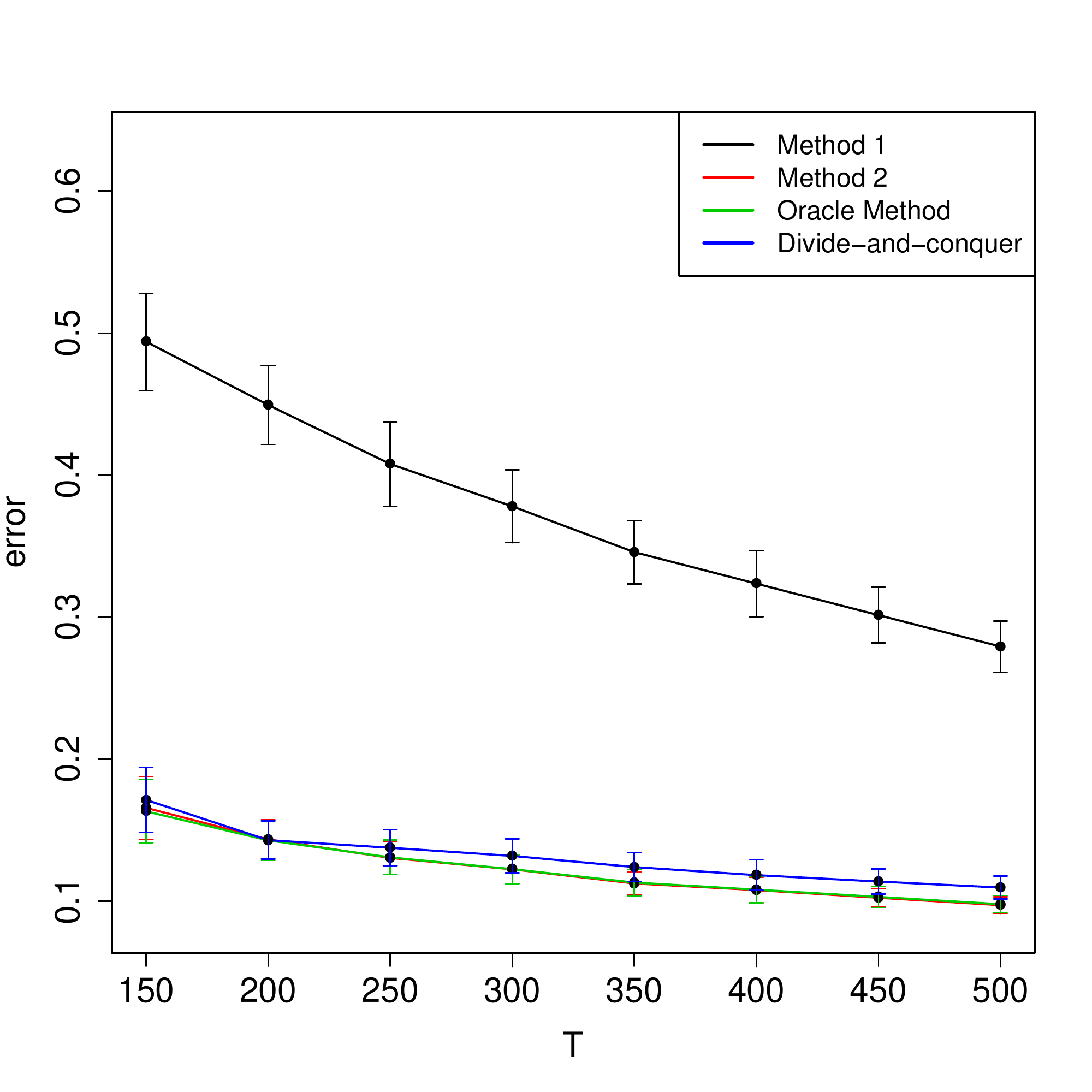}
    \caption{$\snorm{\hat{\Sigmab}_S-\Sigmab_S}$}
  \end{subfigure}%
  \begin{subfigure}[b]{0.5\textwidth}
    \centering
    \includegraphics[width=1.0\textwidth]{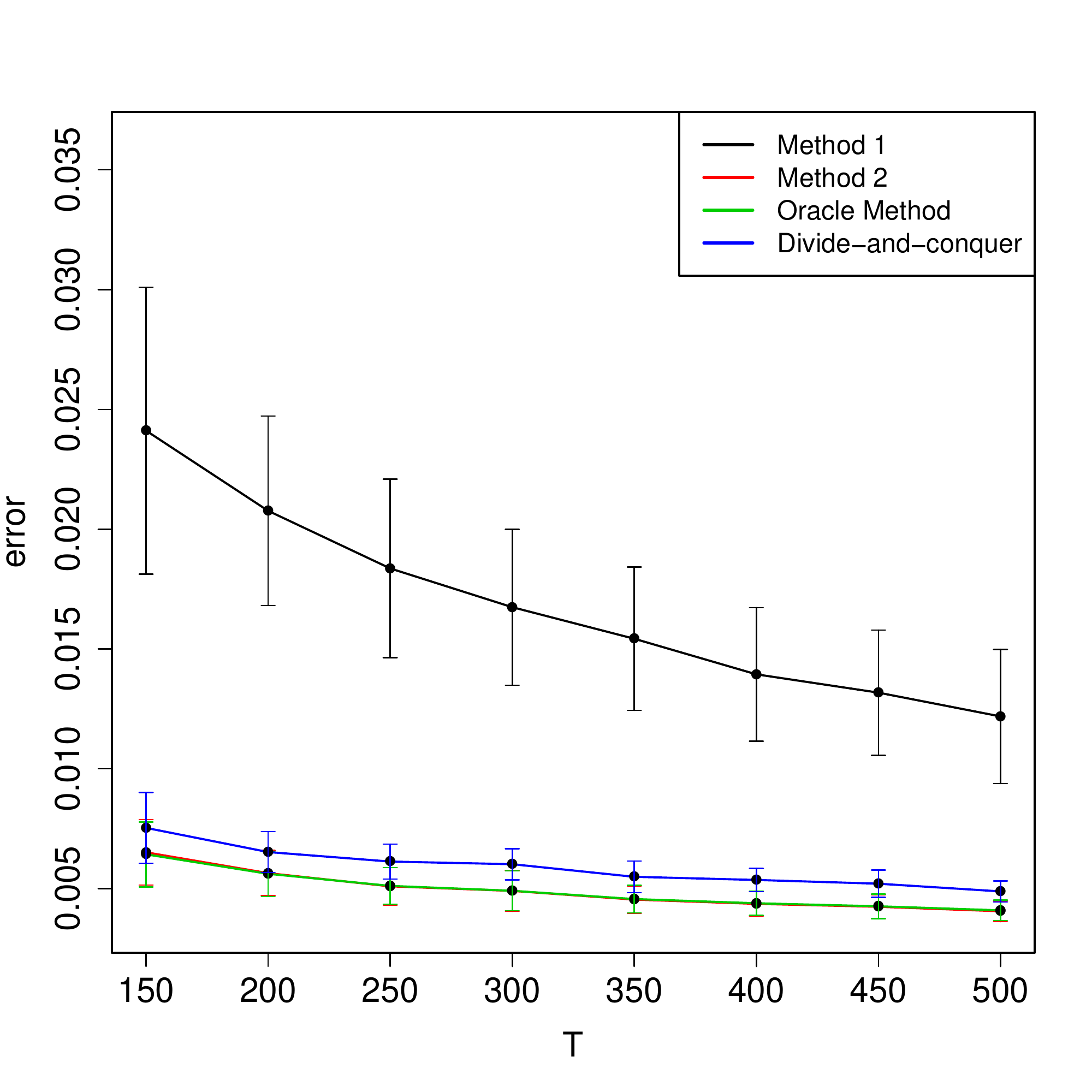}
    \caption{$\norm{\hat{\Sigmab}^{-1}_S-\Sigmab_S^{-1}}$}
  \end{subfigure}
  \begin{subfigure}[b]{0.5\textwidth}
    \centering
    \includegraphics[width=1.0\textwidth]{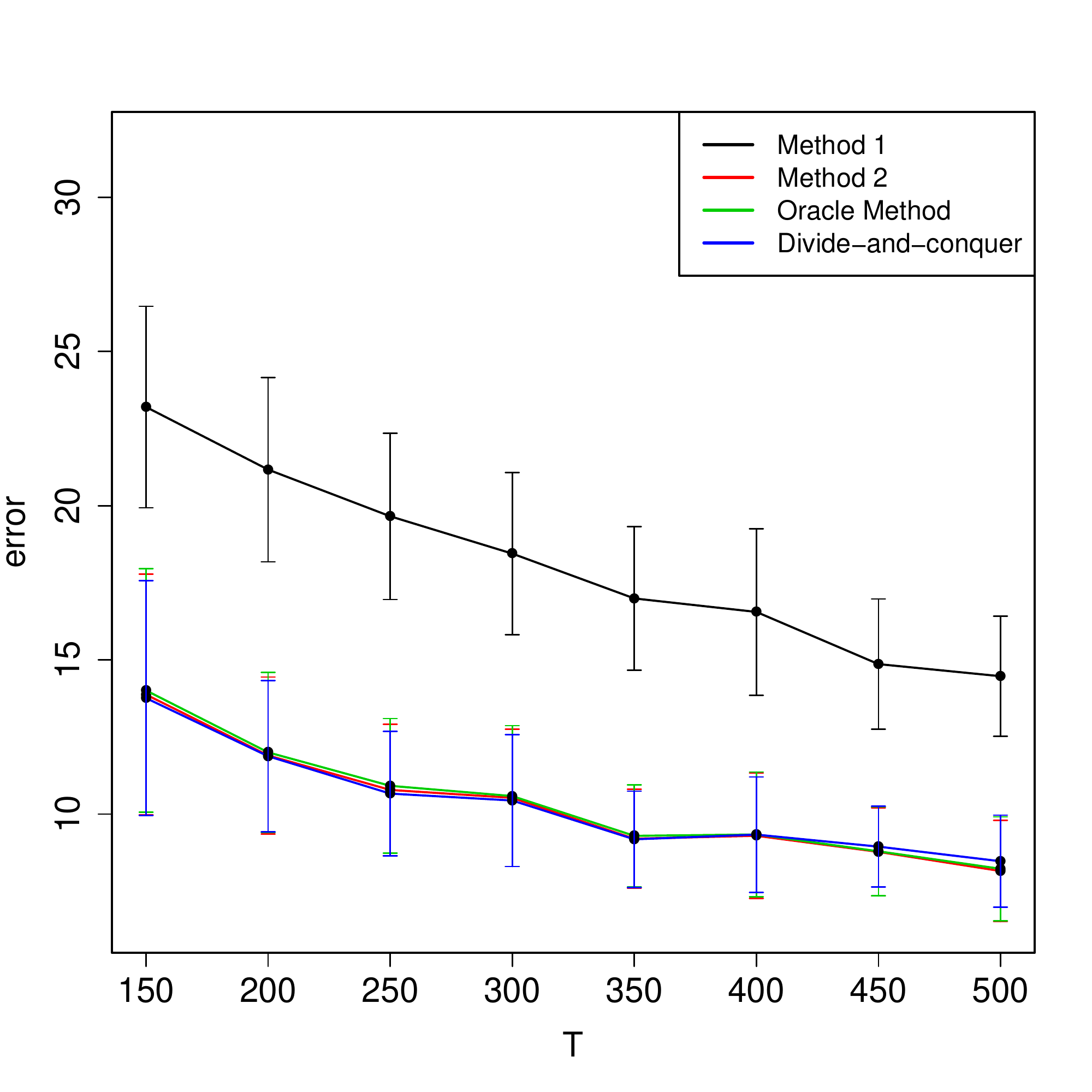}
    \caption{$\maxnorm{\hat{\Sigmab}_S-\Sigmab_S}$}
  \end{subfigure}%
  \begin{subfigure}[b]{0.5\textwidth}
    \centering
    \includegraphics[width=1.0\textwidth]{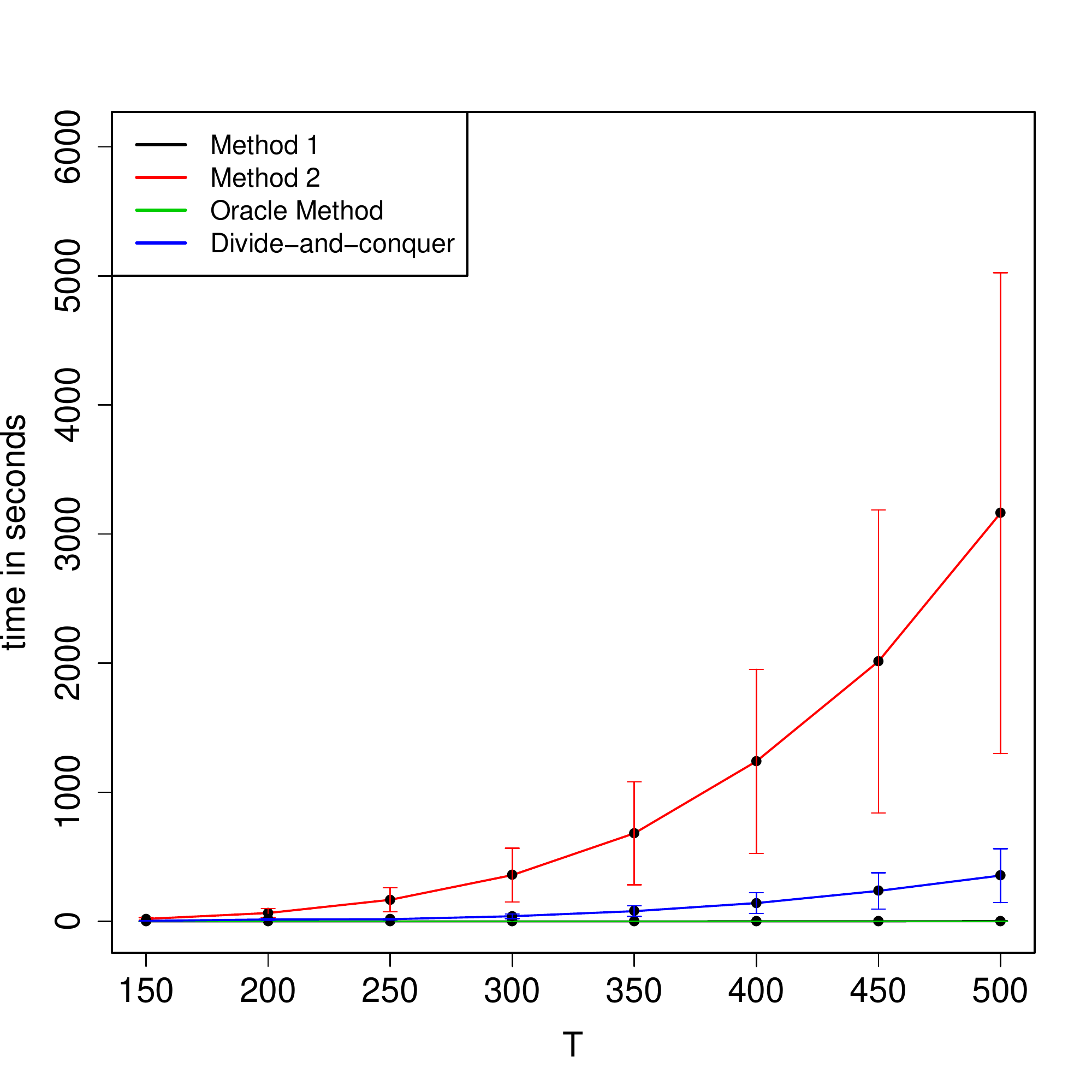}
    \caption{Computational time}
  \end{subfigure}
  \caption{Estimation error by four methods and their computational time: the
    dotted lines represent the means over 100 simulations and the segments
    represent the corresponding standard deviations.}
  \label{fig:1}
\end{figure}

For further comparison with the divide-and-conquer method, we vary $T$ from 50 to
500 and set $(s,p,M)$ as $s=\lfloor T^{0.6}\rfloor$, $p=\lfloor T^{1.4}\rfloor$
and $M=\lfloor T^{0.2}\rfloor$. Figure \ref{fig:1} shows the estimation errors of
the four methods together with the corresponding computational time. Again, when
$p$ is large, Method 2 performs as well as the Oracle Method, both of which
greatly outperform Method 1. However, its computation becomes much slower in this
case. In contrast, the divide-and-conquer method is much faster, while
maintaining comparable performance as Method 2. In the extreme case that $p$ is
around 6000 ($T=500$), the divide-and-conquer method can boost the speed by 9
fold for Method 2.

\begin{table}[hbtp]
  \begin{center} \scriptsize
    \begin{tabular}{c|c|c|c|c|c|c}
      \hline\hline
      $(s, p)$ & \multicolumn{3}{|c|}{$(50,1000)$} & \multicolumn{3}{|c}{$(50, 2000)$}  \\
      \hline
      Method &M1&M2&ORA &M1&M2&ORA \\
      \hline
      $||\hat{\Sigmab}_S - \Sigmab_S||_{\Sigmab_S}$ & 0.271(0.014) & 0.205(0.013) & 0.204(0.013) & 0.270(0.014) & 0.201(0.013) & 0.200(0.013)  \\
      \hline
      $||{\hat{\Sigmab}_S}^{-1} - \Sigmab_S^{-1}||$ & 0.016(0.003) & 0.009(0.002) & 0.009(0.002) & 0.017(0.003) & 0.009(0.002) & 0.009(0.002)  \\
      \hline
      $||\hat{\Sigmab}_S - \Sigmab_S||_{\text{max}}$ & 18.828(3.072) & 17.460(3.237) & 17.457(3.261) & 18.076(2.697) & 16.631(2.949) & 16.623(2.950)  \\
      \hline\hline
      $\max_{t\leq T} ||\hat{\f_t} - \H\f_t ||$ & 1.811(0.195) & 0.445(0.046) & NA & 1.870(0.236) & 0.331(0.025) & NA   \\
      \hline
      $\max_{i\leq s} ||\hat{\b}_i - \H\b_i ||$ & 8.064(0.694) & 4.100(0.330) & 3.858(0.274) & 8.150(0.682) & 3.932(0.292) & 3.805(0.297)    \\
      \hline
      $\max_{i\leq s, t \leq T} ||\hat{\b}_i' \hat{\f_t} - \b_i' \f_t ||$ & 11.375(1.262) & 5.519(0.813) & 5.268(0.843) & 11.466(1.353) & 5.253(0.776) & 5.113(0.739)  \\
      \hline\hline
    \end{tabular}
    \caption{Comparison of three methods when $s$ is much smaller than $p$ (T =
      200). M1, M2 and ORA stand for Method 1, 2 and Oracle method,
      respectively.}
    \label{table:1}
  \end{center}
\end{table}

\begin{table}[hbtp]
  \begin{center} \scriptsize
    \begin{tabular}{c|c|c|c|c|c|c}
      \hline\hline
      $(s, p)$                                                            & \multicolumn{3}{|c|}{$(50,1000)$} & \multicolumn{3}{|c}{$(50,2000)$}  \\
      \hline
      Method                                                              & M1                                & M2     & ORA    & M1      & M2     & ORA \\
      \hline
      $||\hat{\Sigmab}_S - \Sigmab_S||_{\Sigmab_S}$                       & 0.186(0.009) & 0.132(0.007) & 0.131(0.007)  & 0.186(0.009) & 0.131(0.008) & 0.130(0.008)  \\
      \hline
      $||{\hat{\Sigmab}_S}^{-1} - \Sigmab_S^{-1}||$                       & 0.011(0.002) & 0.004(0.001) & 0.004(0.001)  & 0.011(0.002) & 0.004(0.001) & 0.004(0.001)  \\
      \hline
      $||\hat{\Sigmab}_S - \Sigmab_S||_{\text{max}}$                      & 14.054(1.945) & 11.922(2.245) & 11.891(2.262) & 14.180(2.154) & 11.901(2.603) & 11.900(2.604)  \\
      \hline\hline
      $\max_{t\leq T} ||\hat{\f_t} - \H\f_t ||$                           & 1.839(0.193) & 0.417(0.036)   & NA  & 1.843(0.198) & 0.305(0.026)  & NA  \\
      \hline
      $\max_{i\leq s} ||\hat{\b}_i - \H\b_i ||$                           & 6.960(0.584) & 2.830(0.200) & 2.692(0.198)  & 7.024(0.605) & 2.761(0.188) & 2.692(0.194)    \\
      \hline
      $\max_{i\leq s, t \leq T} ||\hat{\b}_i' \hat{\f_t} - \b_i' \f_t ||$ & 11.871(1.540) & 4.138(0.510) & 3.824(0.501)  & 11.457(1.569) & 4.088(0.516) & 3.889(0.542)  \\
      \hline \hline
    \end{tabular}
    \caption{Comparison of three methods when $s$ is much smaller than $p$ (T =
      400). M1, M2 and ORA stand for Method 1, 2 and Oracle method,
      respectively.}
    \label{table:2}
  \end{center}
\end{table}

\begin{table}
  \begin{center} \scriptsize
    \begin{tabular}{c|c|c|c|c|c|c}
      \hline\hline
      $(s, p)$ & \multicolumn{3}{|c|}{$(800, 1000)$} & \multicolumn{3}{|c}{$(800, 2000)$}  \\
      \hline
      Method &M1&M2&ORA &M1&M2&ORA \\
      \hline
      $||\hat{\Sigmab}_S - \Sigmab_S||_{\Sigmab_S}$ & 0.440(0.006) & 0.439(0.006) & 0.435(0.006)  & 0.439(0.006) & 0.436(0.006) & 0.435(0.006)   \\
      \hline
      $||{\hat{\Sigmab}_S}^{-1} - \Sigmab_S^{-1}||$ & 0.062(0.009) & 0.062(0.009) & 0.062(0.009) & 0.061(0.009) & 0.061(0.009) & 0.062(0.012)  \\
      \hline
      $||\hat{\Sigmab}_S - \Sigmab_S||_{\text{max}}$ & 24.565(2.626) & 24.562(2.609) & 24.567(2.599) & 24.511(2.883) & 24.543(2.847) & 24.536(2.851)   \\
      \hline\hline
      $\max_{t\leq T} ||\hat{\f_t} - \H\f_t ||$ & 0.488(0.047) & 0.447(0.040)  & NA & 0.478(0.049) & 0.337(0.038) & NA   \\
      \hline
      $\max_{i\leq s} ||\hat{\b}_i - \H\b_i ||$ & 15.550(0.488) & 15.370(0.462) & 14.418(0.271) & 15.595(0.551) & 15.041(0.357) & 14.398(0.243)   \\
      \hline
      $\max_{i\leq s, t \leq T} ||\hat{\b}_i' \hat{\f_t} - \b_i' \f_t ||$ & 6.745(0.611) & 6.680(0.635) & 6.405(0.630) & 6.904(0.734) & 6.697(0.763) & 6.588(0.737)    \\
      \hline\hline
    \end{tabular}
    \caption{Comparison of three methods when $s$ is comparative to $p$ (T =
      200). M1, M2 and ORA stand for Method 1, 2 and Oracle method,
      respectively.}
    \label{table:3}
  \end{center}
\end{table}

\begin{table}[hbtp]
  \begin{center} \scriptsize
    \begin{tabular}{c|c|c|c|c|c|c}
      \hline\hline
      $(s, p)$ & \multicolumn{3}{|c|}{$(800, 1000)$} & \multicolumn{3}{|c}{$(800, 2000)$}  \\
      \hline
      Method &M1&M2&ORA &M1&M2&ORA \\
      \hline
      $||\hat{\Sigmab}_S - \Sigmab_S||_{\Sigmab_S}$ & 0.193(0.004) & 0.192(0.004) & 0.189(0.004) & 0.192(0.004) & 0.190(0.004) & 0.188(0.004)   \\
      \hline
      $||{\hat{\Sigmab}_S}^{-1} - \Sigmab_S^{-1}||$ & 0.008(0.001) & 0.008(0.001) & 0.008(0.001) & 0.008(0.001) & 0.008(0.001) & 0.008(0.001)  \\
      \hline
      $||\hat{\Sigmab}_S - \Sigmab_S||_{\text{max}}$ & 17.062(2.603) & 17.051(2.612) & 17.041(2.621) & 16.919(2.182) & 16.891(2.206) & 16.888(2.209)  \\
      \hline\hline
      $\max_{t\leq T} ||\hat{\f_t} - \H\f_t ||$ & 0.467(0.038) & 0.423(0.036) & NA & 0.466(0.038) & 0.304(0.026) & NA   \\
      \hline
      $\max_{i\leq s} ||\hat{\b}_i - \H\b_i ||$ & 11.009(0.298) & 10.850(0.302) & 10.225(0.205) & 10.934(0.274) & 10.530(0.213) & 10.189(0.172)   \\
      \hline
      $\max_{i\leq s, t \leq T} ||\hat{\b}_i' \hat{\f_t} - \b_i' \f_t ||$ & 5.367(0.577) & 5.276(0.560) & 4.880(0.528) & 5.293(0.411) & 5.024(0.461) & 4.894(0.420)    \\
      \hline\hline
    \end{tabular}
    \caption{Comparison of three methods when $s$ is comparative to $p$ (T =
      400). M1, M2 and ORA stand for Method 1, 2 and Oracle method,
      respectively.}
    \label{table:4}
  \end{center}
\end{table}

\begin{figure}[hbtp]
  \centering
  \includegraphics[width=0.6\textwidth]{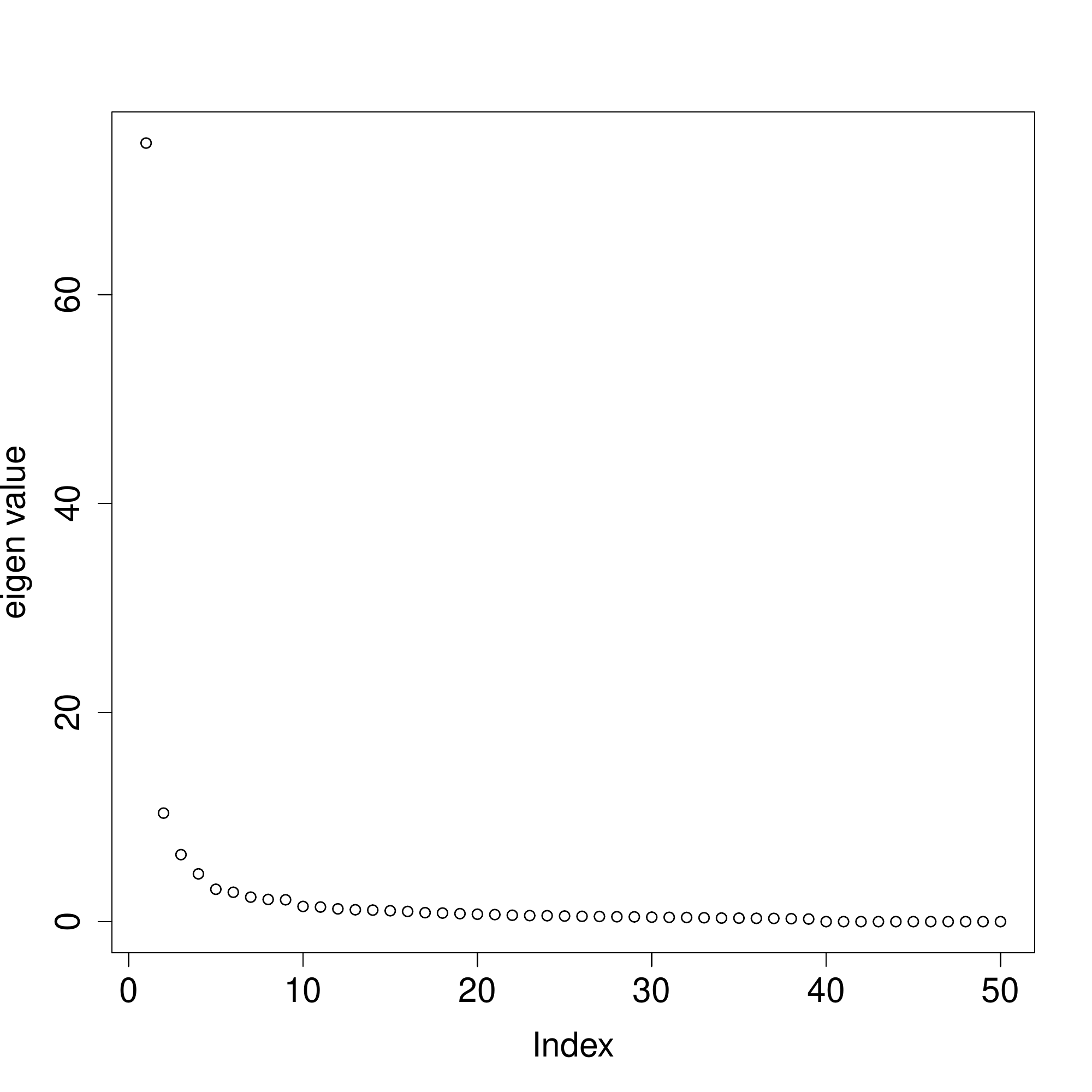}
  \caption{Eigen-values of the sample covariance matrix for GSE22255}
  \label{fig:eigen-value}
\end{figure}
\begin{figure}[hbtp]
  \centering
  \includegraphics[width=0.6\textwidth]{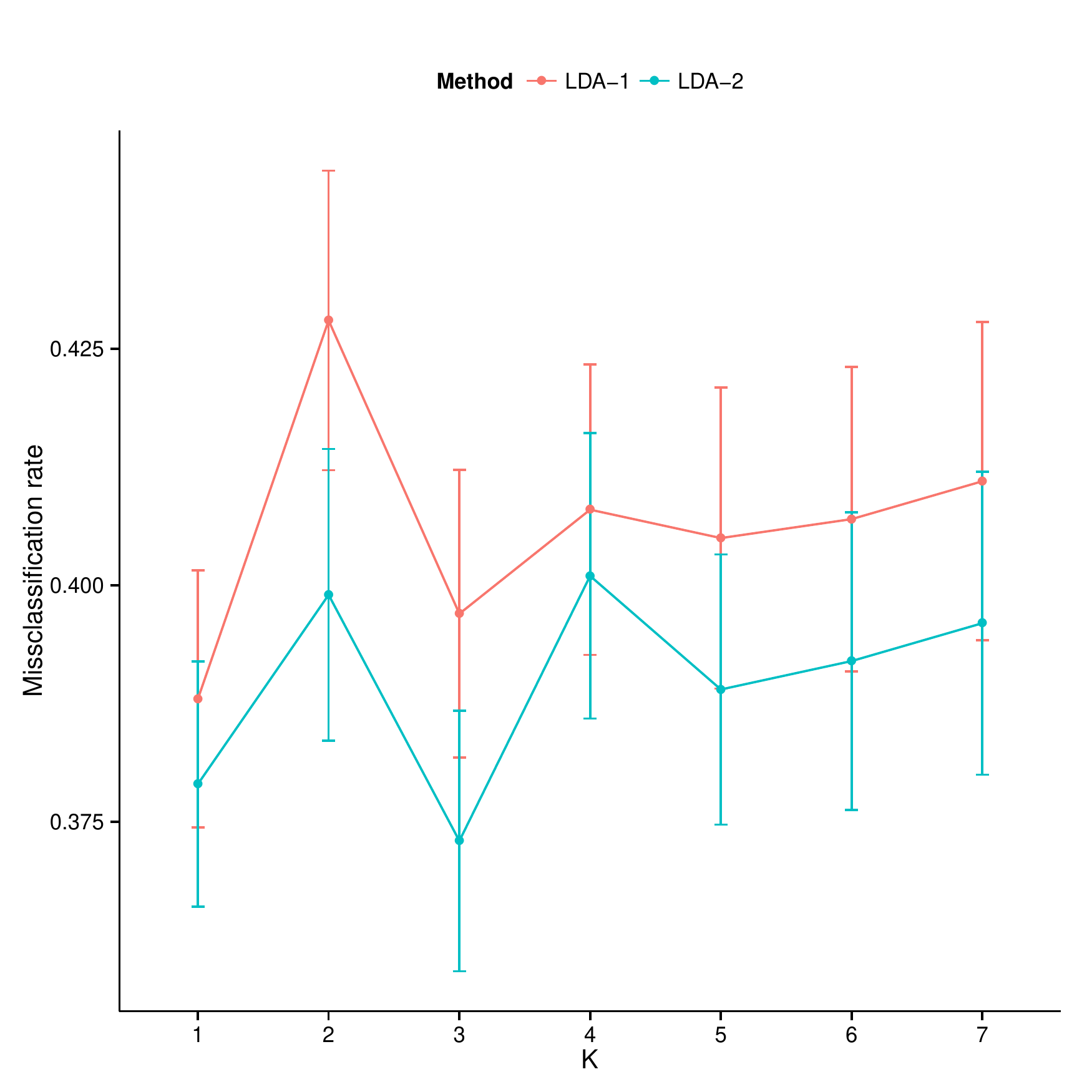}
  \caption{Misclassification rates of LDA-1 and LDA-2 over 100 random splits: the dotted lines represent the
    means over 100 splits and the segments represent the corresponding standard deviations.}
  \label{fig:real-data example}
\end{figure}
\section{Real Data Example}
\label{sec:7}
We use a real data example to illustrate how different utilization of available variables can affect the
inference of the variables of interest. \cite{KGT12} carried out a gene profiling study among 40 Portuguese
and Spanish adults to identify key genetic risk factors for ischemic stroke. Among them, 20 subjects were
patients having ischemic stroke and the others were controls. Their gene profiles were obtained using the
GeneChip Human Genome U133 Plus 2.0 microarray. The data was available at Gene Expression Omnibus with access
name
``GSE22255''. 

To judge how effectively the gene expression can distinguish ischemic stroke and controls, we applied the
Linear Discriminant Analysis (LDA) to this dataset. We randomly chose 10 subjects as the test set and the rest
as the training set. We repeated the random splitting for 100 runs. In each run, we selected the set of
expressed differentially (DE) genes with a threshold of over 1.2-fold change and a Q-value $\leq 0.05$, which
is a commonly used quantity to define DE genes \citep{SJD02}. A LDA rule was then learned from the training
set using the selected genes and further applied to the test set for classifying cases and controls. The LDA
rule classifies a subject as a case if
\begin{equation}
  \label{eq:11}
  \hat{\deltab}'\hat{\Sigmab}^{-1}(\x-\bar{\mub}) \geq 0,
\end{equation}
where $\hat{\deltab}=\hat{\mub}_1-\hat{\mub}_0 \in \mathbb{R}^s$ is the sample mean difference between the
two groups (case - control), $s$ is the number of selected genes, $\hat{\Sigmab}\in \mathbb{R}^{s\times s}$
is an estimator of the true covariance matrix $\Sigmab$ of the selected genes, and
$\bar{\mub}=(\hat{\mub}_1+\hat{\mub}_0)/2$. $\bar{\mub}$, $\hat{\deltab}$ and $\hat{\Sigmab}$ are obtained
from the training set and $\x$ is the gene expression of subjects in the test set.

As $s$ can be larger than the sample size, the traditional LDA where $\hat{\Sigmab}$ is the sample covariance
is no longer applicable. An alternative method to estimate $\Sigmab$ is adopting the factor model. Factor
modeling is widely used in the genomics literature to model the dependencies among genes
\citep{CL12,KS06}. Several factors, like the natural pathway structure \citep{OG00} can be the latent factors
affecting the correlation among genes. A few spiked eigenvalues of the sample covariance in Figure
\ref{fig:eigen-value} also suggest the existence of potential latent factors in this dataset. Again, there are
two ways utilizing the factor model. One way is to use Method 1, where all procedures are done based on the
selected genes only. The resulting rule is referred as ``LDA-1'' in Figure \ref{fig:real-data
  example}. Another way is to use auxiliary data as in Method 2. More specifically, it firstly uses data from
all involved genes and subjects in the training set to estimate the latent factors. These estimated factors
are then applied to the set of selected genes, where their loadings and idiosyncratic matrix estimators are
obtained. Combing them together produces the covariance matrix estimator, which is still an $s\times s$
matrix. The resulting rule is referred as ``LDA-2'' in Figure \ref{fig:real-data example}. Recall that the
only difference between the two rules is that they use different covariance estimators.


Figure \ref{fig:real-data example} plots the average misclassification rates on the test set against the
number of factors for the 100 random splits. It is clearly seen that LDA-2 gives better
misclassification rates than LDA-1, which is solely due to a different estimation of the covariance
matrix. The results lend further support to our claim that using more data is beneficial.

\newpage

\begin{center}
  {\Large\bf Appendix}
\end{center}
\setcounter{section}{0}
\renewcommand{\thesection}{A.\arabic{section}}
\setcounter{equation}{0}
\renewcommand{\theequation}{A.\arabic{equation}}

\section{Additional Regularity Conditions}
\label{sec:A1}
(iv) $\{\u_t,\f_t \}_{t\geq 1}$ are i.i.d. sub-Gaussian random variables over $t$.
\\
(v) There exist constants $c_1$ and $c_2$ that
$0<c_1\leq\lambda_{\min}(\Sigmab_u)\leq \lambda_{\max}(\Sigmab_u)\leq
c_2<\infty$,
$\lonenorm{\Sigmab_u}<c_2$ and $\min_{i\leq p, j\leq p} \var(u_{it}u_{jt})>c_1$; \\
(vi) There exists an $M>0$ such that $||\B||_{\max} <M$;\\
(vii) There exists an $M>0$ such that for any $s\leq T$ and $t\leq T$,
$\E|p^{-1/2} (\u_s'\Sigmab_u^{-1}\u_t-\E\u_s'\Sigmab_u^{-1}\u_t)|^4<M$ and $\E
||p^{-1/2}\B'\Sigmab_u^{-1}\u_t||^4<M$;\\
(viii) For each $t\leq T$, $\E \norm{(pT)^{-1/2} \sum_{s=1}^T\f_s(\u_s'\Sigmab_u^{-1}\u_t-\E(\u_s'\Sigmab_u^{-1}\u_t))}^2=O(1)$;\\
(ix) For each $i \leq p$, $\E \norm{(pT)^{-1/2}\sum_{t=1}^T
  \sum_{j=1}^p\d_j(u_{jt}u_{it}-\E u_{jt}u_{it})}=O(1)$, where $\d_j$ is the $j$th column of $\B'\Sigmab_u^{-1}$;\\
(x) For each $i\leq K$,
$\E \norm{(pT)^{-1/2} \sum_{t=1}^T \sum_{j=1}^N\d_ju_{jt}f_{it}}=O(1)$.

Condition (iv) is a standard assumption in order to establish the exponential type of concentration inequality
for the elements in $\u_t$ and $\f_t$. Condition (v) requires $\Sigmab_u$ to be well-conditioned. In
particular, we need a lower bound on the eigen-values of $\Sigmab_u$. This assumption guarantees that
$\tilde{\Sigmab}_u$ is asymptotically non-singular so that $\tilde{\Sigmab}_u^{-1}$ will not perform badly in
the weighted least-squares problem described in (6). These conditions were also assumed in
\cite{FLM13}. Conditions (vii)-(x) are some moment conditions needed to establish the central limit theorem
for the WPC estimator $\hat{\f}_t$. They are standard in the factor model literature, e.g. \cite{SW02} and
\cite{B03}.

\section{Proofs of Results in Sections \ref{sec:2} and \ref{sec:4}}
\label{sec:A2}
\begin{proof}[\textbf{Proof of Proposition \ref{pro:1}.}]
  Let $\g_1=\nabla_{\thetab_S} \log h(\y_S-\thetab_S,\y_{\Sc}-\thetab_{\Sc})$ and
  $\g_2=\nabla_{\thetab_S}\log h_S(\y_S-\thetab_S)$, where $h_S$ is the marginal density of $\y_S$. Firstly, we show that
  $\g_2=\E(\g_1|\y_S)$. In fact, for any bounded function $\varphi(\y_S)$, by Fubini Theorem and condition (\ref{eq:3}),
  \begin{align*}
    \E(\g_1\varphi(\y_S))
    &=-\iint (\nabla_{\y_S} \log
      h(\y_S-\thetab_S,\y_{\Sc}-\thetab_{\Sc}))h(\y_S-\thetab_S,\y_{\Sc}-\thetab_{\Sc}) \varphi(\y_S)\dd\y_{S}\dd\y_{\Sc}\\
    &=-\iint (\nabla_{\y_S} h(\y_S-\thetab_S,\y_{\Sc}-\thetab_{\Sc}))\varphi(\y_S)\dd\y_{S}\dd\y_{\Sc}\\
    &=-\int \left(\nabla_{\y_S} \int h(\y_S-\thetab_S,\y_{\Sc}-\thetab_{\Sc})\dd\y_{\Sc}\right) \varphi(\y_S)\dd\y_{S}\\
    &=-\int \nabla_{\y_S} h_S(\y_S-\thetab_S)\varphi(\y_S) \dd \y_S\\
    &=\int \left(\nabla_{\y_S} \log h_S(\y_S-\thetab_S) \right)h_S(\y_S-\thetab_S) \varphi(\y_S) \dd\y_S\\
    &=\E(\g_2\varphi(\y_S)).
  \end{align*}
  Then, by definition, $\g_2=\E(\g_1|\y_S)$. Therefore,
  \begin{align*}
    \{I_p(\thetab)\}_S
    &=\E(\g_1\g_1')=\E[(\g_2+\g_1-\g_2)(\g_2+\g_1-\g_2)']\\
    &=\E[\g_2\g_2']+\E[\g_2(\g_1-\g_2)']+\E[(\g_1-\g_2)\g_2']+\E[(\g_1-\g_2)(\g_1-\g_2)']\\
    &=I_S(\thetab_S)+\E[(\g_1-\g_2)(\g_1-\g_2)']\\
    &\succeq I_S(\thetab_S),
  \end{align*}
  where the last equality follows from
  $\E[\g_2(\g_1-\g_2)']=\E[\E[\g_2(\g_1-\g_2)'|\y_S]]=0$, since $\g_2=\E(\g_1|\y_S)$.
\end{proof}

\begin{proof}[\textbf{Proof of Example 2.}]
  Without loss of generality, we assume $\thetab=\bdm{0}$ so that the density of $\y$ is proportional
  to $g(\y'\Omegab\y)$, where $\Omegab=\Sigmab^{-1}$. Then,
  \begin{align*}
    \left|\nabla_{\y_S} h(\y_S,\y_{S^c}) \right|
    &= 2\left|g'(\y'\Omegab\y) (\Omegab\y)_S \right|
      \leq 2\left|g'(\y'\Omegab\y) \right| \left|\Omegab_{S}\y_S+\Omegab_{S,S^c}\y_{S^c} \right|\\
    &\leq 2c \left|\Omegab_{S}\y_S+\Omegab_{S,S^c}\y_{S^c} \right| g(\y'\Omegab\y).
  \end{align*}
  Note that
  \begin{align*}
    \int \left( \int  \left|\Omegab_{S}\y_S+\Omegab_{S,S^c}\y_{S^c} \right| g(\y'\Omegab\y)
    d\y_{S^c}\right)d\y_S
    &\propto\E \left(
      \left|\Omegab_{S}\y_S+\Omegab_{S,S^c}\y_{S^c} \right| \right) \\
    &\leq  \E \left(\left|\Omegab_{S}\y_S\right|+\left|\Omegab_{S,S^c}\y_{S^c} \right| \right)\\
    &<\infty
  \end{align*}
  Therefore for a.e. any $\y_S$, $\int \left|\Omegab_{S}\y_S+\Omegab_{S,S^c}\y_{S^c} \right| g(\y'\Omegab\y)$ is
  integrable. By Example 1.8 of \cite{SJ03}, differentiation and integration are interchangeable, hence
  (\ref{eq:3}) holds.
\end{proof}

\begin{proof}[\textbf{Proof of Proposition \ref{pro:2}.}]
  For simplicity, let $\Omegab=I_p(\thetab)$ and partition it as
  \begin{equation*}
    \Omegab=
    \begin{pmatrix}
      \Omegab_{S} & \Omegab_{S,S^c} \\ \Omegab_{S^c,S} & \Omegab_{\Sc}
    \end{pmatrix}.
  \end{equation*}
  Then, the Fisher information $I(\f)$ of $\f$ contained in all data is given by
  \begin{align}
    I(\f)= \B' \Omegab\B
    =\B_S'\Omegab_{S}\B_S+\B_{S^c}'\Omegab_{S^c,S}\B_S+\B_S'\Omegab_{S,S^c}\B_{S^c}+
    \B_{S^c}'\Omegab_{S^c}\B_{S^c}.\label{eq:A1}
  \end{align}
  If $\Omegab_{S,\Sc}=\bdm{0}$, we have
  \begin{align*}
    I(\f) &=\B_S'\Omegab_{S}\B_S+\B_{S^c}'\Omegab_{S^c}\B_{S^c}
            =\B_S'\{I_p(\thetab) \}_S\B_S+\B_{S^c}'\Omegab_{S^c}\B_{S^c}\\
          &\succeq\B_S'I_S(\thetab_S)\B_S+\B_{S^c}'\Omegab_{S^c}\B_{S^c}
            \succeq \B_S'I_S(\thetab_S)\B_S=I_S(\f),
  \end{align*}
  where the first inequality follows from Proposition \ref{pro:1} and the last inequality
  follows from that $\B_{S^c}'\Omegab_{S^c}\B_{S^c}$ is positive semi-definite. This
  completes the proof.
\end{proof}

\begin{proof}[\textbf{Proof of Proposition \ref{pro:3}.}]
  For any general $\Q\in \mathbb{R}^{L\times R}$, $\B_L\in \mathbb{R}^{L\times K}$, and $\B_R\in \mathbb{R}^{R\times K}$,
  we have
  \begin{equation*}
    \E(\B_L'\Q\B_R)= \E \left[\sum_{l=1}^L \sum_{r=1}^R q_{l,r} \b_{L,l} \b_{R,r}'\right].
  \end{equation*}
  where $q_{l,r}$ is the $(l,r)$-th element of $\Q$, $\b'_{L,l}$ is the $l$th row of $\B_L$ and $\b'_{R,r}$ is the $r$th
  row of $\B_R$. Therefore,
  \begin{equation*}
    \E(\B_{S^c}'\Omegab_{S^c,S}\B_S) = \E \left[\sum_{l\in S^C} \sum_{r\in S} \omega_{l,r}
      \b_{S^c,l}\b_{S,r}' \right],
  \end{equation*}
  where $\omega_{l,r}$ is the $(l,r)$-th element of $\Omegab$. By the i.i.d assumption, for $l \in S^C$ and $r\in S$,
  $\E(\b_{S^c,l}\b_{S,r}')=\E(\b_{S^c,l})\E(\b'_{S,r})=\bdm{0}$. Hence,
  $\E(\B_{S^c}'\Omegab_{S^c,S}\B_S)=\bdm{0}$. Similarly, it can be shown that
  $\E(\B_S'\Omegab_{S,S^c}\B_{S^c})=\bdm{0}$. By Proposition \ref{pro:1}, $\B_S'\Omegab_{S}\B_S \succeq \I_S(\f)$, which
  implies that $\E(\B_S'\Omegab_{S}\B_S) \succeq \E(\I_S(\f))$.
  \begin{equation*}
    \E(\B_{S^c}'\Omegab_{S^c}\B_{S^c})= \E \left[\sum_{l\in S^c} \sum_{r\in S^c} \omega_{l,r} \b_{L,l} \b_{R,r}' \right] =
    \E \left[\sum_{l\in S^c} \omega_{l,l} \b_{L,l}\b_{L,l}' \right]= \mathrm{tr}(\Omegab_{S^c})\E(\b\b') \succeq \bdm{0}.
  \end{equation*}
  Using (\ref{eq:A1}) and the above results, we have $\E[I(\f)] \succeq \E[I_S(\f)]$.
\end{proof}

\begin{proof}[\textbf{Proof of Lemma \ref{lem:1}.}]
  Since we assume all conditions hold for both $s$ and $p$, we prove the result for $p$, i.e.
  $\max_{t\leq T} \norm{\hat{\f}^{(2)}_t-\H_2 \f_t} =O_P\left(T^{-1/2}+T^{1/4}/p^{-1/2}\right)$. The
  result for $s$ can be proved similarly. For simplicity, we write $\hat{\f}_t^{(2)}$ as $\hat{\f}_t$ and $\H_2$ as
  $\H$.

  By (A.1) of \cite{BW13}, $\hat{\f}_t-\H\f_t$ has the following expansion,
  \begin{equation*}
    \hat{\f}_t-\H\f_t=\hat{\V}^{-1} \left(\frac{1}{T} \sum_{i=1}^T \hat{\f}_i\u_i'\wt\u_t/p + \frac{1}{T} \sum_{i=1}^T
      \hat{\f}_i\hat{\eta}_{it} + \frac{1}{T} \sum_{i=1}^T \hat{\f}_i\hat{\theta}_{it} \right),
  \end{equation*}
  where $\hat{\eta}_{it}=\f_i'\B'\wt\u_t/p$,
  $\hat{\theta}_{it}=\f_t'\B'\wt\u_i/p$, and $\hat{\V}$ is the diagonal matrix of
  the $K$ largest eigenvalues of $\Y'\tilde{\Sigmab}_u^{-1}\Y/T$. Let $\eta_{it}=\f_i'\B'\w\u_t/p$ and
  ${\theta}_{it}=\f_t'\B'\w\u_i/p$. Then, we have
  \begin{align}
    \norm{\hat{\f}_t-\H\f_t}
    &\leq \norm{\hat{\V}^{-1}} \left(\bignorm{\frac{1}{T} \sum_{i=1}^T \hat{\f}_i\u_i'(\wt-\w)\u_t/p}
      +\bignorm{\frac{1}{T} \sum_{i=1}^T \hat{\f}_i(\u_i'\w\u_t-\E \u_i'\w\u_t)/p} \right. \nonumber \\
    &\left. \hspace{3ex} +\bignorm{\frac{1}{T} \sum_{i=1}^T\hat{\f}_i\E(\u_i'\w\u_t)/p} + \bignorm{\frac{1}{T} \sum_{i=1}^T
      \hat{\f}_i(\hat{\eta}_{it}-\eta_{it})}+\bignorm{\frac{1}{T}\sum_{i=1}^T \hat{\f}_i\eta_{it}} \right.\nonumber \\
    &\left. \hspace{3ex} +\bignorm{\frac{1}{T} \sum_{i=1}^T
      \hat{\f}_i(\hat{\theta}_{it}-\theta_{it})}+ \bignorm{\frac{1}{T}\sum_{i=1}^T \hat{\f}_i\theta_{it}}  \right). \label{eq:A2}
  \end{align}
  Denote the $j$th summand inside the parenthesis as $G_{jt}$.

  By Lemma A.2 of \cite{BW13}, $\norm{\hat{\V}^{-1}}=O_P\left(1 \right)$. By
  Lemma A.6(iv) of \cite{BW13},
  \[\max_{t\leq T} G_{1t}=O_P\left(\norm{\wt-\w}\left\{\norm{\wt-\w}+1/\sqrt{p}+\sqrt{(\log p)/T} \right\} \right).\] By Proposition
  4.1 of \cite{BW13},
  \begin{equation}
    \label{eq:A3}
    \norm{\wt-\w}=o_P\left(\min \left\{T^{-1/4},p^{-1/4}, \sqrt{T/(p\log p)} \right\} \right),
  \end{equation}
  therefore, $\norm{\wt-\w}\left(\norm{\wt-\w}+1/\sqrt{p}+\sqrt{(\log p)/T}\right)=o(T^{-1/2}+p^{-1/2})$. Hence,
  \begin{equation*}
    \max_{t\leq T} G_{1t}= o_P\left(T^{-1/2}+p^{-1/2} \right).
  \end{equation*}
  By Lemma A.8(ii) of \cite{BW13},
  $\max_{t\leq T} G_{2t}=O_P\left(T^{1/4}p^{-1/2} \right)$. By Lemma A.10(i) of
  \cite{BW13}, $\max_{t\leq T} G_{3t}=O_P\left(T^{-1/2} \right)$. By Lemma A.6(vi) of
  \cite{BW13},
  \begin{equation*}
    \max_{t\leq T} G_{4t}=O_P\left(\norm{\wt-\w}\left\{\norm{\wt-\w}+1/\sqrt{p}+1/\sqrt{T} \right\} \right)+o_P\left(1/\sqrt{p} \right)=
    o_P\left(1/\sqrt{p} \right).
  \end{equation*}
  By Lemma A.8(iii) of \cite{BW13}, $\max_{t \leq T}G_{5t}=O_P\left(T^{1/4}p^{-1/2} \right)$. By Lemma A.6(v) of
  \cite{BW13} and (\ref{eq:A3}),
  \begin{equation*}
    \max_{t\leq T} G_{6t}=O_P\left(\norm{\wt-\w}\left\{\norm{\wt-\w}+1/\sqrt{p} + \sqrt{(\log p)/T}\right\} \right)  =o_P\left(1/\sqrt{p} \right).
  \end{equation*}
  By Lemma A.6(iii) of \cite{BW13} and (\ref{eq:A3}),
  \begin{equation*}
    \max_{t \leq T} G_{7t}=O_P\left( \norm{\wt-\w}/\sqrt{p}+1/p+1/\sqrt{pT} \right)=o_P\left( 1/\sqrt{p} \right).
  \end{equation*}
  Then, by (\ref{eq:A2}), we have
  \begin{equation*}
    \max_{t\leq T} \norm{\hat{\f}_t-\H\f_t}=O_P\left( \frac{1}{\sqrt{T}}+ \frac{T^{1/4}}{\sqrt{p}} \right).
  \end{equation*}
\end{proof}

\begin{proof}[\textbf{Proof of Lemma \ref{lem:2}.}]
  For Method 1, we have the following decomposition
  \begin{equation*}
    \hat{\b}^{(1)}_i-\H_1\b_i=\underbrace{\frac{1}{T}\sum_{t=1}^T\H_1\f_tu_{it}}_{I_1}+
    \underbrace{\frac{1}{T} \sum_{t=1}^Ty_{it}(\hat{\f}^{(1)}_t-\H_1\f_t)}_{I_2} +
    \underbrace{\H_1(\frac{1}{T}\sum_{t=1}^T\f_t\f_t'-\I_K)\b_i}_{I_3},
  \end{equation*}
  where $\b_i$ is the true factor loading of the $i$th subject as defined in (\ref{eq:1}).

  For $I_1$, we have
  \begin{equation*}
    \max_{i\leq s} \bignorm{\frac{1}{T} \sum_{t=1}^T \H_1\f_tu_{it}} \leq \norm{\H_1} \max_{i\leq
      s} \sqrt{\sum_{k=1}^K \Big(\frac{1}{T} \sum_{t=1}^T f_{kt}u_{it}\Big)^2}.
  \end{equation*}
  It follows from Lemma C.3(iii) of \cite{FLM13} that,
  $\max_{i\leq s} \sqrt{\sum_{k=1}^K (\frac{1}{T} \sum_{t=1}^T
    f_{kt}u_{it})^2}\allowbreak=O_P\left(\sqrt{(\log s)/T} \right)$.
  From Lemma \ref{lem:A2}, $\norm{\H_1}=O_P\left(1 \right)$, therefore
  $I_1=O_P\left(\sqrt{(\log s)/T} \right)$.

  As for $I_2$, by conditions (v) and (vi),
  \begin{align*}
    \max_{i\leq s} \E y_{it}^2=\max_{i\leq s} \{\E(\b_i'\f_t)^2+\E u_{it}^2 \}\leq\max_{i\leq s}
    \norm{\b_i}^2+\max_{i \leq s}
    \var(u_{it})=O(1).
  \end{align*}
  By condition (iv), $y_{it}^2$ is sub-exponential, therefore by the union bound and sub-exponential tail bound,
  $ \max_{i\leq s} \left|\frac{1}{T} \sum_{t=1}^T y_{it}^2-\E y_{it}^2 \right| =O_P\left(\sqrt{(\log s)/T} \right)$.
  Then,
  \begin{equation}
    \label{eq:A4}
    \max_{i\leq s} \frac{1}{T} \sum_{t=1}^T y_{it}^2\leq \max_{i\leq s} \left|\frac{1}{T} \sum_{t=1}^T y_{it}^2-\E
      y_{it}^2 \right| + \max_{i\leq s} \E y_{it}^2=O_P\left(1 \right).
  \end{equation}
  By Cauchy-Schwartz inequality,
  \begin{align*}
    \max_{i\leq s} \bignorm{\frac{1}{T} \sum_{t=1}^T y_{it}(\hat{\f}^{(1)}_t-\H_1\f_t)}
    &\leq \max_{i\leq s} \left(\frac{1}{T}
      \sum_{t=1}^T y_{it}^2 \cdot \frac{1}{T} \sum_{t=1}^T
      \norm{\hat{\f}^{(1)}_t-\H_1\f_t}^2
      \right)^{1/2}\\
    &=O_P\left(\left(\frac{1}{T} \sum_{t=1}^T
      \norm{\hat{\f}^{(1)}_t-\H_1\f_t}^2\right)^{1/2}\right)\\
    &=O_P\left(\frac{1}{\sqrt{T}}+\frac{1}{\sqrt{s}} \right),
  \end{align*}
  where the last equality follows from Lemma \ref{lem:A5}. So, $I_2=O_P\left(1/\sqrt{T}+1/\sqrt{s} \right)$.

  Finally, it follows from Lemma C.3(i) of \cite{FLM13} that
  $\norm{\frac{1}{T} \sum_{t=1}^T\f_t\f_t'-\I_K}=O_P\left(T^{-1/2} \right)$. This together with
  $\norm{\H_1}=O_P\left(1 \right)$ and condition (vi) show that $I_3=O_P\left(T^{-1/2} \right)$. Hence,
  \begin{equation*}
    \max_{i\leq s} \norm{\hat{\b}_i^{(1)}-\H_1\b_i}=O_P\left(\frac{1}{\sqrt{s}}+\sqrt{\frac{\log s}{T}}\right).
  \end{equation*}

  Using the same arguments and the results of $\hat{\f}_t^{(2)}$ in Lemma \ref{lem:1}, we can show that
  \begin{equation*}
    \max_{i\leq s} \norm{\hat{\b}_i^{(2)}-\H_2\b_i}=O_P\left(\frac{1}{\sqrt{p}}+\sqrt{\frac{\log s}{T}}\right).\\
  \end{equation*}
  When the common factor $\f_t$ is known, for the oracle estimator of the loading matrix, we have
  \begin{align*}
    \max_{i\leq s}\norm{\hat{\b}^o_i-\b_i}
    &\leq\max_{i\leq s}\bignorm{\frac{1}{T} \sum_{t=1}^T \f_tu_{it}}+\bignorm{\frac{1}{T}
      \sum_{t=1}^T \f_t\f_t'-\I_K}\max_{i\leq s}
      \norm{\b_i}\\&=O_P\left(\sqrt{\frac{\log
                     s}{T}}+\frac{1}{\sqrt{T}}  \right)\\
    &= O_P\left(\sqrt{\frac{\log s}{T}} \right).
  \end{align*}
\end{proof}

\begin{proof}[\textbf{Proof of Lemma \ref{lem:3}.}]
  By Theorem A.1 of \cite{FLM13} (cited as Lemma \ref{lem:A7} in Appendix), it suffices to show
  \begin{equation*}
    \max_{i\leq s} \frac{1}{T} \sum_{t=1}^T (u_{it}
    -\hat{u}^{(1)}_{it})^2=O_P\left(\frac{1}{s}+\frac{\log s}{T}
    \right) \quad \text{ and } \quad \max_{i,t} |u_{it}-\hat{u}^{(1)}_{it}|=o_P\left(1 \right).
  \end{equation*}
  For Method 1, we have
  \begin{equation*}
    u_{it}-\hat{u}^{(1)}_{it}=\b_i'\H_1'(\hat{\f}^{(1)}_t-\H_1\f_t) +
    \{(\hat{\b}^{(1)}_i)'-\b_i'\H_1 \}\hat{\f}^{(1)}_t+\b_i'(\H_1'\H_1-\I_K)\f_t
  \end{equation*}
  Using $(a+b+c)^2\leq4a^2+4b^2+4c^2$, we have
  \begin{align*}
    \max_{i\leq s} \frac{1}{T} \sum_{t=1}^T (u_{it}-\hat{u}^{(1)}_{it})^2
    &\leq 4
      \max_{i\leq s} \norm{\H_1\b_{i}}^2
      \frac{1}{T} \sum_{t=1}^T \norm{\hat{\f}_t^{(1)}-\H_1 \f_t}^2 \\
    &+ 4\max_{i\leq s}
      \norm{\hat{\b}^{(1)}_i-\H_1\b_i}^2 \frac{1}{T}
      \sum_{t=1}^T \norm{\hat{\f}^{(1)}_t}^2 \\
    &+4 \fnorm{\H_1'\H_1-\I_K}^2\max_{i\leq s} \norm{\b_i}^2\frac{1}{T} \sum_{t=1}^T
      \norm{\f_t}^2.
  \end{align*}
  Since,
  $\max_i \norm{\H_1\b_i}\leq \norm{\H_1}\max_i \norm{\b_i}=O_P\left(1 \right)$,
  $\frac{1}{T}\sum_{t=1}^T \norm{\hat{\f}_t^{(1)}}^2=O_P\left(1 \right)$, and
  $\frac{1}{T} \sum_{t=1}^T \norm{\f_t}^2\allowbreak=O_P\left(1 \right)$, it follows from
  Lemma \ref{lem:1}, \ref{lem:2}, \ref{lem:A3} and \ref{lem:A5} that
  \begin{equation}
    \label{eq:A5}
    \max_{i\leq s} \frac{1}{T} \sum_{t=1}^T (u_{it}
    -\hat{u}^{(1)}_{it})^2=O_P\left(\frac{1}{s}+\frac{\log s}{T}
    \right).
  \end{equation}
  On the other hand, by Lemma \ref{lem:A1},
  \begin{equation*}
    \max_{i,t} |u_{it}-\hat{u}^{(1)}_{it}|=\max_{i,t}
    |(\hat{\b}^{(1)}_i)'\hat{\f}^{(1)}_i-\b'_i\f_t|= O_P\left(
      (\log T)^{1/2} \sqrt{\frac{\log s}{T}}+
      \frac{T^{1/4}}{\sqrt{s}} \right)=o(1).
  \end{equation*}
  Then, the result follows from Theorem A.1 of \cite{FLM13}.

  In analogous, a similar result can be proved for Method 2. For the oracle estimator,
  $\hat{u}^o_{it}=y_{it}-(\hat{\b}_i^o)'\f_t$. Therefore,
  \begin{gather*}
    \max_{i\leq s} \frac{1}{T} \sum_{t=1}^T (u_{it}-\hat{u}^{o}_{it})^2 \leq \max_{i\leq s}
    \norm{\hat{\b}_i^o-\b_i}^2 \frac{1}{T} \sum_{t=1}^T \norm{\f_t}^2=O_P\left(\max_{i\leq s}
      \norm{\hat{\b}_i^o-\b_i}^2 \right)= O_P\left(\frac{\log s}{T} \right).\\
    \max_{i,t} |u_{it}-\hat{u}_{it}^{o}|=\max_{i,t} |(\hat{\b}_i^o)'\f_t-\b'_i\f_t| =
    O_P\left( (\log T)^{1/2}\sqrt{\frac{\log s}{T}} \right)=o_P(1).
  \end{gather*}
  It then follows from Theorem A.1 of \cite{FLM13} that
  \begin{equation*}
    \norm{\hat{\Sigmab}_{u,S}^{o}-\Sigmab_{u,S}}=O_P\left(m_s\sqrt{\frac{\log s}{T}} \right)=\norm{(\hat{\Sigmab}_{u,S}^{o})^{-1}-\Sigmab_{u,S}^{-1}}.
  \end{equation*}
\end{proof}

\begin{proof}[\textbf{Proof of Theorem \ref{thm:1}.}]
  (1) For Method 1, $\hat{\Sigmab}^{(1)}_S=\hat{\B}_1\hat{\B}_1'+\hat{\Sigmab}^{(1)}_{u,S}$. Therefore,
  \begin{align*}
    \snorm{\hat{\Sigmab}^{(1)}_S-\Sigmab_S}^2
    &\leq
      2 \left(\snorm{\hat{\B}_1\hat{\B}_1'-\B_S\B_S'}^2+
      \snorm{\hat{\Sigmab}^{(1)}_{u,S}-\Sigmab_{u,S}}^2 \right)\\
    &\leq 2
      \left(\snorm{\B_{S}(\H_1'\H_1-\I_K)\B_S'}^2+2\snorm{\B_S\H'_{1}\C_1'}^2 +\snorm{\C_1\C_1'}^2
      \right. \\
    &\hspace{3ex}\left.+
      \snorm{\hat{\Sigmab}^{(1)}_{u,S}-\Sigmab_{u,S}}^2 \right),
  \end{align*}
  where $\C_1=\hat{\B}_1-\B_S\H_1'$. Then, it follows from Lemmas \ref{lem:A4} that
  \begin{align*}
    \snorm{\hat{\Sigmab}_S^{(1)}-\Sigmab_S}^2=O_P\left(\frac{1}{sT}+\frac{1}{s^2}+w_1^2+sw_1^4+m_s^2w_1^2
    \right)
    =O_P\left(sw_1^4+m_s^2w_1^2 \right).
  \end{align*}
  Similarly,
  $ \snorm{\hat{\Sigmab}_S^{(2)}-\Sigmab_S}^{2}=O_P\left(sw_2^4+m_s^2w_2^2 \right)$.

  In the oracle case, we have
  \begin{align*}
    \snorm{\hat{\Sigmab}_S^o-\Sigmab_S}^2
    &\leq
      2 \left(\snorm{\hat{\B}_o\hat{\B}_o'-\B_S\B_S'}^2+
      \snorm{\hat{\Sigmab}^{o}_{u,S}-\Sigmab_{u,S}}^2 \right)\\
    &\leq
      2 \Big(\underbrace{\snorm{(\hat{\B}_o-\B_S)(\hat{\B}_o-\B_S)'}^2}_{I_1} +
      2\underbrace{\snorm{(\hat{\B}_o-\B_S)\B_S'}^2}_{I_2} +
      \underbrace{\snorm{\hat{\Sigmab}^{o}_{u,S}-\Sigmab_{u,S}}^2}_{I_3} \Big).
  \end{align*}
  Since all eigenvalues of $\Sigmab_S$ are bounded away from zero, for any matrix $\A\in \mathbb{R}^{s\times s}$,
  $\snorm{\A}^2=s^{-1}\fnorm{\Sigmab^{-1/2}\A\Sigmab^{-1/2}}^2=O_P\left(s^{-1} \fnorm{\A}^2 \right)$. Then, by Lemma
  \ref{lem:2}, we have
  \begin{equation*}
    I_1=O_P\left(s^{-1} \fnorm{\hat{\B}_o-\B_S}^4 \right) = O_P\left(s w_o^4\right),
  \end{equation*}
  where the last equality follows that $\fnorm{\hat{\B}_o-\B_S}^2\leq s(\max_{i\leq s}
  \norm{\hat{\b}_i^o-\b_i})^2=O_P\left(sw_o^2 \right)$. For $I_2$, we have
  \begin{align*}
    I_2&=s^{-1} \tr((\hat{\B}_o-\B_S)'\Sigmab_S^{-1}(\hat{\B}_o-\B_S) \B_S'\Sigmab_S^{-1}\B_S) \\
       &\leq s^{-1} \norm{\Sigmab_S^{-1}} \fnorm{\hat{\B}_o-\B_S}^2 \norm{\B_S'\Sigmab_S^{-1}\B_S}\\
       &=O_P\left(w_o^2 \right).
  \end{align*}
  For $I_3$, Lemma \ref{lem:3} implies that
  \begin{equation*}
    I_3=O_P\left(s^{-1} \fnorm{\hat{\Sigmab}_{u,S}^o-\Sigmab_{u,S}}^2\right)
    =O_P\left(\norm{\hat{\Sigmab}_{u,S}^o-\Sigmab_{u,S}}^2 \right)=O_P\left(m_s^2 w_o^2 \right).
  \end{equation*}
  Therefore, $\snorm{\hat{\Sigmab}_{u,S}^o-\Sigmab_{u,S}}^2=O_P\left(sw_o^4+m_s^2w_o^2\right)$.

  (2) For Method 1,
  \begin{equation*}
    \maxnorm{\hat{\Sigmab}_S^{(1)}-\Sigmab_S} \leq \underbrace{\maxnorm{\hat{\B}_1\hat{\B}_1'-\B_S\B_S'}}_{I_1}+
    \underbrace{\maxnorm{\hat{\Sigmab}^{(1)}_{u,S}-\Sigmab_{u,S}}}_{I_2}.
  \end{equation*}
  For $I_1$, we have
  \begin{align*}
    I_1&=\max_{ij} |(\hat{\b}_i^{(1)})'\hat{\b}_j^{(1)}-\b_i'\b_j|\\
       &\leq \max_{ij} \left(|(\hat{\b}_i^{(1)}-\H_1\b_i)'(\hat{\b}_j^{(1)}-\H_1\b_j)|+2|\b_i'\H_1'(\hat{\b}_j^{(1)}-\H_1\b_j)|+
         |\b_i'(\H_1\H_1'-\I_K)\b_j| \right)\\
       &\leq \big(\max_i \norm{\hat{\b}_i^{(1)}-\H_1\b_i}\big)^2+2\max_{ij}
         \norm{\hat{\b}_i^{(1)}-\H_1\b_i}\norm{\H_1\b_j}+\norm{\H_1\H_1'-\I_K} \big(\max_i \norm{\b_i}\big)^2\\
       &=O_P\left( w_1 \right),
  \end{align*}
  where the last identity follows from Lemmas \ref{lem:2} and \ref{lem:A3}.

  For $I_2$, let $\sigma_{u,ij}$ be the $(i,j)$-th entry of
  $\Sigmab_{u,S}$ and $\hat{\sigma}_{u,ij}=\frac{1}{T}
  \sum_{t=1}^T \hat{u}_{it}\hat{u}_{jt}$, where $\hat{u}_{it}$ are
  the estimator of $u_{it}$ from Method 1 as described in Section
  \ref{sec:4}. Then,
  \begin{align*}
    &\hspace{3ex}\max_{ij} |\hat{\sigma}_{u,ij}-\sigma_{u,ij}|\\
    &=\max_{ij}\Big|\frac{1}{T} \sum_{t=1}^T
      (\hat{u}_{it}\hat{u}_{jt}-u_{it}u_{jt})\Big|+
      \max_{ij} \Big|\frac{1}{T} \sum_{i=1}^T u_{it}u_{jt}-\E(u_{it}u_{jt})\Big|\\
    &\leq \max_{ij} \Big|\frac{1}{T} \sum_{t=1}^T (\hat{u}_{it}-u_{it})(\hat{u}_{jt}-u_{jt})\Big| + 2\max_{ij} \Big|\frac{1}{T}
      \sum_{t=1}^T (\hat{u}_{it}-u_{it})u_{jt}\Big|
      + \max_{ij} \Big|\frac{1}{T} \sum_{i=1}^T
      u_{it}u_{jt}-\E(u_{it}u_{jt})\Big|\\
    &\leq \max_{ij} \left(\frac{1}{T} \sum_{t=1}^T(\hat{u}_{it}-u_{it})^2 \right)^{1/2} \left(\frac{1}{T} \sum_{t=1}^T
      (\hat{u}_{jt}-u_{jt})^2 \right)^{1/2}
      + 2\max_{ij}\left(\frac{1}{T} \sum_{t=1}^T(\hat{u}_{it}-u_{it})^2 \right)^{1/2} \left(\frac{1}{T}
      \sum_{t=1}^T u_{jt}^2\right)^{1/2} \\
    &\hspace{3ex}+\max_{ij} \Big|\frac{1}{T} \sum_{i=1}^T  u_{it}u_{jt}-\E(u_{it}u_{jt})\Big|\\
    &=O_P\left(w_1^2 \right)+O_P\left(w_1 \right)+O_P\left(\sqrt{(\log s)/T} \right),
  \end{align*}
  where the last equality follows from (\ref{eq:A5}), Lemma C.3 (ii) of \cite{FLM13} and
  $$\max_{j\leq s} \frac{1}{T} \sum_{t=1}^T u_{jt}^2=O_P\left(1 \right)$$ as similarly shown in
  (\ref{eq:A4}). Hence, $\max_{ij}|\hat{\sigma}_{u,ij}-\sigma_{u,ij}|=O_P\left(w_1 \right)$. After the
  thresholding,
  \begin{align*}
    \max_{ij}|s_{ij}(\hat{\sigma}_{u,ij})-\sigma_{u,ij}|
    &\leq
      \max_{ij} |s_{ij}(\hat{\sigma}_{u,ij})-\hat{\sigma}_{u,ij}|+
      |\hat{\sigma}_{u,ij}-\sigma_{u,ij}|\\
    & \leq
      \max_{ij}|s_{ij}(\hat{\sigma}_{u,ij})-\hat{\sigma}_{u,ij}|+O_P\left(w_1 \right)\\
    &=O_P\left(w_1 \right).
  \end{align*}
  where $s_{ij}(\cdot)$ is the hard thresholding at the level defined in step ii. of Method 1. Hence,
  $\maxnorm{\hat{\Sigmab}_{u,S}^{(1)}-\Sigmab_{u,S}}=O_P\left(w_1 \right)$. Similarly,
  $\maxnorm{\hat{\Sigmab}_{u,S}^{(2)}-\Sigmab_{u,S}}=O_P\left(w_2 \right)$.  For the oracle estimator,
  \begin{align*}
    \maxnorm{\hat{\B}_o\hat{\B}_o'-\B\B'}
    &=\max_{ij} \left(|(\hat{\b}_i^o-\b_i)'(\hat{\b}_i-\b_i)|+
      2|(\hat{\b}_i^o-\b_i)'\b_j|\right)\\
    &\leq \Big(\max_i \norm{\hat{\b}_i^o-\H_1\b_i}\Big)^2 + 2\max_{ij} \norm{\hat{\b}_i^o-\b_i} \norm{\b_j} \\
    &=O_P\left(w_o \right),
  \end{align*}
  where the last equality follows from condition (vi) and Lemma \ref{lem:2}. Using similar arguments as in the
  above, $\max_{ij}|\hat{\sigma}_{u,ij}^o-\sigma_{u,ij}|=O_P(w_0)$. Hence,
  $\maxnorm{\hat{\Sigmab}_{u,S}^o-\Sigmab_{u,S}}=O_P\left(w_o \right)$.

  (3) For Method 1, let $\tilde{\Sigmab}_S=\B_S\H_1'\H_1\B_S'+\Sigmab_{u,S}$. We have
  \begin{equation*}
    \norm{(\hat{\Sigmab}_S^{(1)})^{-1}-\Sigmab_S^{-1}} \leq \norm{(\hat{\Sigmab}_S^{(1)})^{-1}-\tilde{\Sigmab}_S^{-1}}+ \norm{\tilde{\Sigmab}_S^{-1}-\Sigmab_S^{-1}}.
  \end{equation*}
  Since $\hat{\Sigmab}_S^{(1)}=\hat{\B}_1\hat{\B}_1'+ \hat{\Sigmab}_{u,S}^{(1)}$, by Sherman-Morrison-Woodbury
  formula,
  \begin{align*}
    \tilde{\Sigmab}_S^{-1}&=\Sigmab_{u,S}^{-1}+\Sigmab_{u,S}^{-1}\B_S\H_1'\G^{-1}\H_1\B_S\Sigmab_{u,S}^{-1}, \\
    (\hat{\Sigmab}^{(1)}_S)^{-1}&=(\hsigus{(1)})^{-1}+(\hsigus{(1)})^{-1}\hat{\B}_1\hat{\G}^{-1}\hat{\B}_1(\hsigus{(1)})^{-1},
  \end{align*}
  where $\G=\GG$ and $\hat{\G}=\GGhat$. Therefore, $\norm{(\hat{\Sigmab}_S^{(1)})^{-1}-\tilde{\Sigmab}_S^{-1}}\leq
  \sum_{i=1}^6 I_i$, where
  \begin{align*}
    I_1&= \norm{(\hsigus{(1)})^{-1}-\Sigmab_{u,S}^{-1}},\\
    I_2&= \norm{\{(\hsigus{(1)})^{-1}-\Sigmab_{u,S}^{-1} \}\hat{\B}_1\hat{\G}^{-1}\hat{\B}_1'(\hsigus{(1)})^{-1}},\\
    I_3&= \norm{\{(\hsigus{(1)})^{-1}-\Sigmab_{u,S}^{-1} \}\hat{\B}_1\hat{\G}^{-1}\hat{\B}_1'\Sigmab_{u,S}^{-1}},\\
    I_4&= \norm{\Sigmab_{u,S}^{-1}(\hat{\B}_1-\B_S\H_1')\hat{\G}^{-1}\hat{\B}_1'\Sigmab_{u,S}^{-1}},\\
    I_5&= \norm{\Sigmab_{u,S}^{-1}(\hat{\B}_1-\B_S\H_1')\hat{\G}^{-1}\H_1\B_S'\Sigmab_{u,S}^{-1}},\\
    I_6&= \norm{\Sigmab_{u,S}^{-1}\B_S\H_1'\{\hat{\G}^{-1}-\G^{-1}\}\H_1\B_S'\Sigmab_{u,S}^{-1}}.
  \end{align*}
  From Lemma \ref{lem:3}, $I_1=O_P\left(m_sw_1 \right)$. For $I_2$, we have
  \begin{equation*}
    I_2\leq \norm{(\hsigus{(1)})^{-1}-\Sigmab_{u,S}^{-1}}
    \norm{\hat{\B}_1\hat{\G}^{-1}\hat{\B}'_1}
    \norm{(\hsigus{(1)})^{-1}}.
  \end{equation*}
  By Lemma \ref{lem:3} and condition (v), $\norm{(\hsigus{(1)})^{-1}}=O_P\left(1 \right)$. Lemma \ref{lem:A6}(ii) implies
  that $\norm{\hat{\G}^{-1}}=O_P\left(s^{-1} \right)$. Therefore,
  $\norm{\hat{\B}_1\hat{\G}^{-1}\hat{\B}_1'}=O_P\left(1 \right)$ and $I_2=O_P\left(m_sw_1 \right)$. Similarly,
  $I_3=O_P\left(m_sw_1 \right)$. For $I_4$, condition (v) implies that $\norm{\Sigmab_{u,S}^{-1}}=O(1)$. Next,
  $\norm{(\hat{\B}_1-\B_S\H_1')\hat{\G}^{-1}\hat{\B}_1'}$ is bounded by
  \begin{equation*}
    \norm{(\hat{\B}_1-\B_S\H_1')\hat{\G}^{-1}\hat{\B}_1'}\leq
    \norm{(\hat{\B}_1-\B_S\H_1')\hat{\G}^{-1}(\hat{\B}_1-\B_S\H_1')'}^{1/2} \norm{\hat{\B}_1\hat{\G}^{-1}\hat{\B}_1'}^{1/2}.
  \end{equation*}
  Since $\norm{\hat{\G}^{-1}}=O_P(s^{-1})$ by Lemma \ref{lem:A6}(ii) and
  $\fnorm{\hat{\B}_1-\B_S\H_1'}^2=O_P\left(sw_1^2 \right)$ by Lemma \ref{lem:A4}(i), we have
  $\norm{(\hat{\B}_1-\B_S\H_1')\hat{\G}^{-1}(\hat{\B}_1-\B_S\H_1')'}= O_P\left(w_1^2 \right)$. This together with
  $\norm{\hat{\B}_1\hat{\G}^{-1}\hat{\B}_1'}\allowbreak=O_P\left(1 \right)$ imply that $I_4=O_P\left(w_1 \right)$. Similarly,
  $I_5=O_P\left(w_1 \right)$. For $I_6$, we have
  \begin{equation*}
    I_6\leq \norm{\Sigmab_{u,S}^{-1}\B_S\H_1'\H_1\B_S'\Sigmab_{u,S}^{-1}} \norm{\hat{\G}^{-1}-\G^{-1}}.
  \end{equation*}
  Condition (ii), (v) and $\norm{\H_1}=O_P\left(1 \right)$ imply that
  $\norm{\Sigmab_{u,S}^{-1}\B_S\H_1'\H_1\B_S'\Sigmab_{u,S}^{-1}}=O_P\left(s \right)$. Next, we bound
  $\norm{\hat{\G}^{-1}-\G^{-1}}$. Note that,
  \begin{align*}
    \norm{\hat{\G}^{-1}-\G^{-1}}
    &=\norm{\G^{-1}(\hat{\G}-\G)\hat{\G}^{-1}}
      =O_P\left(s^{-2}\norm{\hat{\B}_1'(\hsigus{(1)})^{-1}\hat{\B}_1-(\B_S\H_1')'\Sigmab_{u,S}^{-1}\B_S\H_1'}
      \right) \\
    &=O_P\left(s^{-1}m_sw_1 \right),
  \end{align*}
  because by Lemma \ref{lem:A6} (i) and (ii), $\norm{\G^{-1}}=O\left(s^{-1} \right)$,
  $\norm{\hat{\G}^{-1}}=O_P\left(s^{-1} \right)$, and
  \begin{align}
    \label{eq:A6}
    &\hspace{3ex}\norm{\hat{\B}_1'(\hsigus{(1)})^{-1}\hat{\B}_1-(\B_S\H_1')'\Sigmab_{u,S}^{-1}\B_S\H_1'} \nonumber\\
    &\leq \norm{(\hat{\B}_1-\B_S\H_1')' (\hsigus{(1)})^{-1}(\hat{\B}_1-\B_S\H_1')}+2
      \norm{(\hat{\B}_1-\B_S\H_1')(\hsigus{(1)})^{-1} \B_S\H_1'}\nonumber\\
    &\hspace{3ex}+ \norm{(\B_S\H_1')'\{ (\hsigus{(1)})^{-1}-\Sigmab_{u,S}^{-1}\}\B_S\H_1'} \nonumber\\
    &=O_P\left(sw_1^2 \right)+O_P\left(sw_1 \right)+O_P\left(sm_sw_1 \right) \nonumber \\
    &=O_P\left(sm_sw_1 \right).
  \end{align}
  Therefore, $I_6=O_P\left(m_sw_1 \right)$. Summing the six terms, we have
  $\norm{(\hat{\Sigmab}_{u,S}^{(1)})^{-1}-\tilde{\Sigmab}_S^{-1}}=O_P\left(m_sw_1 \right)$. Next, we bound
  $\norm{\tilde{\Sigmab}_S^{-1}-\Sigmab_S^{-1}}$.

  By using Sherman-Morrison-Woodbury formula again,
  \begin{align*}
    \norm{\tilde{\Sigmab}_S^{-1}-\Sigmab_S^{-1}}
    &= \bignorm{\Sigmab_{u,S}^{-1}\B_S
      \big\{[(\H_1'\H_1)^{-1}+\B_S'\Sigmab_{u,S}^{-1}\B_S]^{-1} - [\I_K+\B_S'\Sigmab_{u,S}^{-1}\B_S]^{-1}\big \}\B'_S\Sigmab_{u,S}^{-1}}\\
    &=O(s) \bignorm{[(\H_1'\H_1)^{-1}+\B_S'\Sigmab_{u,S}^{-1}\B_S]^{-1}-[\I_K+\B_S'\Sigmab_{u,S}^{-1}\B_S]^{-1}}\\
    &=O_P\left(s^{-1} \right) \norm{(\H_1'\H_1)^{-1}-\I_K}\\
    &=o_P\left(m_sw_1 \right).
  \end{align*}
  Therefore, $\norm{(\hat{\Sigmab}_{u,S}^{(1)})^{-1}-\Sigmab_S^{-1}}=O_P\left(m_sw_1 \right)$. A similar result can be shown that
  $\norm{(\hat{\Sigmab}_{u,S}^{(2)})^{-1}-\Sigmab_S^{-1}}=O_P\left(m_sw_2 \right)$.

  For the oracle estimator, by Sherman-Morrison-Woodbury formula, $\norm{(\hat{\Sigmab}_S^o)^{-1}-\Sigmab_S^{-1}}\leq
  \sum_{i=1}^6I_i$, where
  \begin{align*}
    I_1&= \norm{(\hsigus{o})^{-1}-\Sigmab_{u,S}^{-1}},\\
    I_2&= \norm{\{(\hsigus{o})^{-1}-\Sigmab_{u,S}^{-1}\}\hat{\B}_o\hat{\J}^{-1}\hat{\B}_o'(\hsigus{o})^{-1}},\\
    I_3&= \norm{\{(\hsigus{o})^{-1}-\Sigmab_{u,S}^{-1} \}\hat{\B}_o\hat{\J}^{-1}\hat{\B}_o'\Sigmab_{u,S}^{-1}},\\
    I_4&= \norm{\Sigmab_{u,S}^{-1}(\hat{\B}_o-\B_S)\hat{\J}^{-1}\hat{\B}_o'\Sigmab_{u,S}^{-1}},\\
    I_5&= \norm{\Sigmab_{u,S}^{-1}(\hat{\B}_o-\B_S)\hat{\J}^{-1}\B_S'\Sigmab_{u,S}^{-1}},\\
    I_6&= \norm{\Sigmab_{u,S}^{-1}\B_S\{\hat{\J}^{-1}-\J^{-1}\}\B_S'\Sigmab_{u,S}^{-1}},
  \end{align*}
  that $\hat{\J}=\I_K+\hat{\B}_o'(\hat{\Sigmab}_{u,S}^o)^{-1}\hat{\B}_o$ and $\J=\I_K+\B_S'\Sigmab_{u,S}^{-1}\B_S$.

  By Lemma \ref{lem:3}, $I_1=O_P\left(m_sw_o \right)$. For $I_2$, Lemma \ref{lem:A6}(ii) implies that
  $\norm{\hat{\J}^{-1}}=O_P\left(s^{-1} \right)$. This together with condition (ii) imply that
  $\norm{\hat{\B}_o\hat{\J}^{-1}\hat{\B}'_o} =O_P\left(1 \right)$. Moreover, it follows from Lemma \ref{lem:3} and
  condition (v) that $\norm{(\hsigus{o})^{-1}}=O_P\left(1 \right)$. Therefore,
  \begin{align*}
    I_2&\leq \norm{(\hsigus{o})^{-1}-\Sigmab_{u,S}^{-1}} \norm{\hat{\B}_o\hat{\J}^{-1}\hat{\B}'_o}
         \norm{(\hsigus{o})^{-1}}=O_P\left(m_sw_o \right).
  \end{align*}
  Similarly, $I_3=O_P\left(m_sw_o \right)$. For $I_4$, we have $I_4\leq
  \norm{(\hat{\B}_o-\B_S)\hat{\J}^{-1}\B_S'}\norm{\Sigmab_{u,S}^{-1}}^2$. We bound
  $\norm{(\hat{\B}_o-\B_S)\hat{\J}^{-1}\B_S'}$ by
  \begin{equation*}
    \norm{(\hat{\B}_o-\B_S)\hat{\J}^{-1}\B_S'} \leq \norm{(\hat{\B}_o-\B_S)\hat{\J}^{-1}(\hat{\B}_o-\B_S)'}^{1/2}
    \norm{{\B}_S\hat{\J}^{-1}{\B}_S'}^{1/2}.
  \end{equation*}
  Since
  $\norm{(\hat{\B}_o-\B_S)(\hat{\B}_o-\B_S)'}\leq \fnorm{\hat{\B}_o-\B_S}^2 \leq s (\max_s
  \norm{\hat{\b}_i^o-\b_i})^2=O_P\left(sw_o^2 \right)$.
  This together with $\norm{\hat{\J}^{-1}}=O_P\left(s^{-1} \right)$ and
  $\norm{\hat{\B}_o\hat{\J}^{-1}\hat{\B}_o}=O_P\left(1 \right)$ imply that $I_4=O_P\left(w_o \right)$. Similarly,
  $I_5= O_P\left(w_o \right)$. For $I_6$, we have
  $I_6\leq \norm{\hat{\J}^{-1}-\J^{-1}}\norm{\Sigmab_{u,S}^{-1}}^2\norm{\B_S\B_S'}$. By conditions (ii) and (iv), we have
  $\norm{\Sigmab_{u,S}^{-1}}=O(1)$ and $\norm{\B_S\B_S'}=O(s)$. As for $\norm{\hat{\J}^{-1}-\J^{-1}}$, we have
  \begin{equation*}
    \norm{\hat{\J}^{-1}-\J^{-1}}=\norm{\hat{\J}^{-1}(\hat{\J}-\J)\J^{-1}} =
    O_P\left(s^{-2}\norm{\B_S'\Sigmab_{u,S}^{-1}\B_S-\hat{\B}_o'\hat{\Sigmab}_{u,S}^{-1}\hat{\B}_o}
    \right)=O_P\left(s^{-1}m_sw_o \right),
  \end{equation*}
  where the last equation follows from that
  \begin{align*}
    \norm{\hat{\B}_o'\hat{\Sigmab}_{u,S}^{-1}\hat{\B}_o -\B_S'\Sigmab_{u,S}^{-1}\B_S}
    & \leq \norm{(\hat{\B}_o-\B_S)' \hat{\Sigmab}_{u,S}^{-1}(\hat{\B}_o-\B_S)}+ 2
      \norm{(\hat{\B_o}-\B_S)'\hat{\Sigmab}_{u,S}^{-1}\B_S} \\
    &+ \norm{\B_S'\{(\hat{\Sigmab}_{u,S}^o)^{-1}-\Sigmab_{u,S}^{-1} \}\B_S} \\
    &=O_P\left(s w_o^2 \right)+O_P\left( sw_o \right)+O_P\left(sm_sw_o \right)\\
    &=O_P\left(sm_sw_o \right).
  \end{align*}
  Therefore, $I_6=O_P\left(m_sw_o \right)$. After summing up,
  $\norm{(\hat{\Sigmab}_S^o)^{-1}-\Sigmab_S^{-1}}=O_P\left(m_sw_o \right)$.
\end{proof}

\section{Convergence Rates of $\bar{\Sigmab}_S$ in Section \ref{sec:5}}
\label{sec:A3}
Let $\bar{\H}=M^{-1} \sum_{m=1}^M \H_{[m]}$, where
$\H_{[m]}=\hat{\V}_m^{-1}\hat{\F}_m'\F_m\B_m'\tilde{\Sigmab}_{u,m}^{-1} \B_m/T$,
$\hat{\V}_m$ is the diagonal matrix of the $K$ largest eigenvalues of
$\Y_m'\tilde{\Sigmab}_{u,m}^{-1}\Y_m/T$, $\B_m$ and $\F_m$ are the loadings and
the factors in the $m$th group.

According to the proof of Theorem \ref{thm:1}, the key is to show that
$\max_{1\leq t\leq T} \norm{\bar{\f}_{t}-\bar{\H}\f_t}$ has the same rate as
$\max_{1\leq t\leq T} \norm{\hat{\f}^{(2)}_{t}-\H_2\f_t}$ and
$\max_{i\leq s} \norm{\bar{\b}_i-\bar{\H}\b_i}$ has the same rate as
$\max_{1\leq i \leq s} \norm{\hat{\b}_i^{(2)}-\H_2\b_i}$.

To give the rate of $\max_{1\leq t\leq T} \norm{\bar{\f}_{t}-\bar{\H}\f_t}$,
since $M$ is fixed, $p/M$ is in the same order as $p$. Then, it follows from
Lemma \ref{lem:1} that for any $1\leq m\leq M$,
$\max_{1\leq t\leq T}\norm{\hat{\f}_{m,t}-\H_{[m]}\f_t}=O_P(a_{p,T})$, where
$a_{p,T}=T^{-1/2}+T^{1/4}p^{-1/2}$. By definition, there exists a positive
constant $C_{m,\epsilon}$ such that
\begin{equation*}
  P \left(\max_{1\leq t\leq T}\norm{\hat{\f}_{m,t}-\H_{[m]}\f_t}>C_{m,\epsilon}a_{p,T}\right) \leq \epsilon/M.
\end{equation*}
Let $C=\max_{1\leq m \leq M} C_{m,\epsilon}$. We have
\begin{align*}
  P \left(\max_{1\leq t\leq T} \norm{\bar{\f}_t-\bar{\H}\f_t}> C a_{p,T}\right)
  &= P \left(\max_{1\leq
    t\leq T} \bignorm{\frac{1}{M} \sum_{m=1}^M (\hat{\f}_{m,t}-\H_{[m]}\f_t)} > C a_{p,T}\right) \\
  &\leq \sum_{m=1}^M P \left(\max_{1\leq t\leq T} \norm{\hat{\f}_{m,t}-\H_{[m]}\f_t} > C a_{p,T}\right) \\
  &\leq \epsilon.
\end{align*}
By definition, $\max_{1\leq t\leq T}\norm{\bar{\f}_t-\bar{\H}\f_t}= O_P\left(a_{p,T}\right)$, which is the same as
$\max_{1\leq t \leq T} \norm{\hat{\f}_t^{(2)}-\H_2\f_{t}}$ shown in Lemma \ref{lem:1}.

Next, we show that $\max_{i\leq s} \norm{\bar{\b}_i-\bar{\H}\b_i}= O_P(w_2)$.
For any $1\leq m \leq M$, similarly as in Lemma \ref{lem:A2}, we have
$\norm{\H_{[m]}}=O_P\left(1 \right)$. By the same union bound argument, we have
$\norm{\bar{\H}}=O_P\left(1 \right)$. Then, it follows from the same proof of
Lemma \ref{lem:2} that
$\max_{i\leq s} \norm{\bar{\b}_i-\bar{\H}\b_i}=O_P\left(w_2 \right)$.

As $M$ is fixed, the results in Lemma \ref{lem:3} and Theorem \ref{thm:1} for
each individual group hold. Repeatedly using the above union bound argument,
$\bar{\Sigmab}_S$ is shown to have the same convergence rate as
$\hat{\Sigmab}_S^{(2)}$.

\section{Additional Lemmas}
\label{sec:A4}
\setcounter{lem}{0}
\renewcommand{\thelem}{A.\arabic{lem}}
\begin{lem}
  \label{lem:A1}
  Under conditions of Lemma \ref{lem:1}, it holds that
  \begin{align*}
    \max_{i\leq s,t\leq T}
    \norm{(\hat{\b}_i^{(1)})'\hat{\f}_t^{(1)}-\b_i'\f_t}
    &=O_P\left((\log
      T)^{1/2}\sqrt{\frac{\log s}{T}} +\frac{T^{1/4}}{\sqrt{s}} \right)\\
    \max_{i\leq s,t\leq T}
    \norm{(\hat{\b}_i^{(2)})'\hat{\f}_t^{(2)}-\b_i'\f_t}
    &=O_P\left((\log
      T)^{1/2}\sqrt{\frac{\log s}{T} } +\frac{T^{1/4}}{\sqrt{p}}\right)\\
    \max_{i\leq s,t\leq T}
    \norm{(\hat{\b}_i^{o})'\f_t-\b_i'\f_t}
    &=O_P\left( (\log
      T)^{1/2} \sqrt{\frac{\log s}{T}}  \right).
  \end{align*}
\end{lem}
\begin{proof}[\textbf{Proof of Lemma \ref{lem:A1}.}]
  Under condition (i), it follows from the union bound argument that
  $$\max_{t\leq T} \norm{\f_t}=O_P\left(\sqrt{\log T} \right).$$ Then, for Method 1, it
  follows from Lemmas \ref{lem:1}, \ref{lem:2}, \ref{lem:A2}, and condition (vi) that, uniformly in $i$ and
  $t$,
  \begin{align*}
    \norm{(\hat{\b}_i^{(1)})'\hat{\f}^{(1)}_t-\b'_i\f_t}
    &\leq \norm{\hat{\b}^{(1)}_i-\H_1\b_i}
      \norm{\hat{\f}^{(1)}_t-\H_1\f_t} + \norm{\H_1\b_i}
      \norm{\hat{\f}^{(1)}_t-\H_1\f_t}\\
    &\hspace{3ex}+  \norm{\hat{\b}^{(1)}_i-\H_1\b_i} \norm{\H_1\f_t}+\norm{\b_i}\norm{\f_t} \fnorm{\H_1'\H_1-\I_K}\\
    &=O_P\left((\log T)^{1/2} \sqrt{\frac{\log s}{T}} +
      \frac{T^{1/4}}{\sqrt{s}}\right).
  \end{align*}
  For Method 2, similar arguments give
  \begin{equation*}
    \max_{i\leq s,t\leq T}
    \norm{(\hat{\b}_i^{(2)})'\hat{\f}_t^{(2)}-\b_i'\f_t}=O_P\left((\log T)^{1/2} \sqrt{\frac{\log s}{T}} +
      \frac{T^{1/4}}{\sqrt{p}}\right).
  \end{equation*}
  In the oracle setting, where the factors are known, we have
  \begin{align*}
    \max_{i\leq s,t\leq T} \norm{(\hat{\b}_i^{o})'\f_t-\b_i'\f_t}
    &=  \max_{i\leq s,t\leq T}\norm{\hat{\b}^{o}_i-\b_i}
      \norm{\f_t}=O_P\left(\sqrt{\log T}\max_{i\leq s} \norm{\hat{\b}_i^o-\b_i}
      \right)\\&=O_P\left( (\log T)^{1/2} \sqrt{\frac{\log s}{T}}\right).
  \end{align*}
\end{proof}

\begin{lem}
  \label{lem:A2}
  Let $\H_1=\hat{\V}_1^{-1}\hat{\F}^{(1)'}\F\B_S'\tilde{\Sigmab}_{u,S}^{-1}\B_S/T$ and
  $\H_2=\hat{\V}_2^{-1}\hat{\F}^{(2)'}\F\B'\tilde{\Sigmab}_u^{-1}\B/T$, where $\hat{\V}_1$ is the diagonal matrix of
  the largest $K$ eigenvalues of $\Y_S'\tilde{\Sigmab}^{-1}_{u,S}\Y_S/T$ and $\hat{\V}_2$ is the diagonal matrix of the
  largest $K$ eigenvalues of $\Y'\tilde{\Sigmab}^{-1}_u\Y/T$. Under conditions of Lemma \ref{lem:1},
  $\norm{\H_1}=O_P\left(1 \right)$ and $\norm{\H_2}=O_P\left(1 \right)$.
\end{lem}

\begin{proof}[\textbf{Proof of Lemma \ref{lem:A2}.}]
  Since $\Sigmab_{u,S}$ is a submatrix of $\Sigmab_u$, it follows from condition (v) that
  $\lambda_{\min}(\Sigmab_{u,S}^{-1})\geq c_2^{-1}$. By Proposition 4.1 of \cite{BW13},
  $\norm{\tilde{\Sigmab}_{u,S}^{-1}-\Sigmab^{-1}_{u,S}}=o_P\left(1 \right)$. Therefore, with probability tending to 1,
  $\norm{\tilde{\Sigmab}_{u,S}^{-1}}\geq 1/(2c_2)$. Then,
  \begin{equation*}
    T^{-1}\Y_S'\tilde{\Sigmab}_{u,S}^{-1}\Y_S =T^{-1}\Y_S'(\tilde{\Sigmab}_{u,S}^{-1}-(1/2c_2)\I)\Y_S+1/(2c_2T)\Y_S'\Y_S.
  \end{equation*}
  Under the pervasive condition (i), it follows from Lemma C.4 of \cite{FLM13} that the $K$th largest eigenvalue of
  $T^{-1}\Y_S'\Y_S$ is larger than $Ms$. Since $T^{-1}\Y_S'(\tilde{\Sigmab}_{u,S}^{-1}-(1/2c_2)\I)\Y_S$ is semi-positive
  definite, it follows from Weyl's inequality that
  \begin{equation*}
    \lambda_K(T^{-1}\Y_S'\tilde{\Sigmab}_{u,S}^{-1}\Y_S)\geq \lambda_K(1/(2c_2T)\Y_S'\Y_S)\geq
    Ms/(2c_2).
  \end{equation*}
  Hence $||\hat{\V}_1^{-1}||=O_P\left(s^{-1} \right)$. Also,
  $\lambda_{\max}(\norm{\F'\F})=\lambda_{\max}(\norm{\sum_{t=1}^T \f_t\f_t'})=O_P\left(T \right)$.  In
  addition,
  $\lambda_{\max} (\norm{\sum_{t=1}^T \hat{\f}_t^{(1)}(\hat{\f}_t^{(1)})'})=O_P\left(T \right)$, where
  the last equation follows from the constraint in (\ref{eq:6}). Then,
  $\norm{(\hat{\F}^{(1)})'\F}\leq \norm{(\hat{\F}^{(1)})'\hat{\F}^{(1)}}^{1/2} \norm{\F'\F}^{1/2}=O_P\left(T
  \right)$.
  These results together with $\norm{\B_S'\tilde{\Sigmab}_{u,S}^{-1}\B_S}=O(s)$ imply that
  $\norm{\H_1}=O_P\left(1 \right)$. Similarly, $\norm{\H_2}=O_P\left(1 \right)$.
\end{proof}

\begin{lem}
  \label{lem:A3}
  (i) $\fnorm{\H_1\H_1'-\I_K}=O_P\left(\frateone\right)$; (ii) $\fnorm{\H_2\H_2'-\I_K}=O_P\left(\fratetwo \right)$.  (iii)
  $\fnorm{\H_1'\H_1-\I_K}=O_P\left(\frateone\right)$; (iv) $\fnorm{\H_2'\H_2-\I_K}=O_P\left(\fratetwo\right)$.
\end{lem}
\begin{proof}[\textbf{Proof of Lemma \ref{lem:A3}.}]
  Let $\hat{\cov}(\H_1\f_t)=\frac{1}{T} \sum_{t=1}^T (\H_1\f_t)(\H_1\f_t)'$. Then,
  \begin{equation*}
    \fnorm{\H_1\H_1'-\I_K}\leq \underbrace{\fnorm{\H_1\H_1'-\hat{\cov}(\H_1\f_t)}}_{I_1}+
    \underbrace{\fnorm{\hat{\cov}(\H_1\f_t)-\I_K}}_{I_2}.
  \end{equation*}
  For $I_1$, we have $I_1\leq \norm{\H_1}^2 \fnorm{\I_K-\hat{\cov}(\f_t)}$, where
  $\hat{\cov}(\f_t)=\frac{1}{T} \sum_{t=1}^T \f_t\f_t'$. It follows from Lemma C.3(i) of \cite{FLM13} that
  $\fnorm{\I_K-\hat{\cov}(\f_t)}=O_P\left(1/\sqrt{T}\right)$. Then, $I_1=O_P\left(1/\sqrt{T}\right)$, since
  $\norm{\H_1}=O_P\left(1\right)$.  For $I_2$, by the identifiability constraint in (\ref{eq:6}), $\frac{1}{T} \sum_{t=1}^T
  \hat{\f}_t^{(1)}\hat{\f}_t^{(1)'}=\I_K$. Therefore,
  \begin{align*}
    I_2& = \bigfnorm{\frac{1}{T} \sum_{t=1}^T \H_1\f_t(\H_1\f_t)' - \frac{1}{T} \sum_{t=1}^T \hat{\f}_t^{(1)}\hat{\f}_t^{(1)'}} \\
       &\leq \bigfnorm{\frac{1}{T} \sum_{t=1}^T (\H_1\f_t-\hat{\f}_t^{(1)})(\H_1\f_t)'}+ \bigfnorm{\frac{1}{T} \sum_{t=1}^T
         \hat{\f}_t^{(1)}(\hat{\f}_t^{(1)}-\H_1\f_t)'}\\
       &\leq \left(\frac{1}{T} \sum_{t=1}^T \norm{\H_1\f_t-\hat{\f}_t^{(1)}}^2\cdot \frac{1}{T} \sum_{t=1}^T \norm{\H_1\f_t}^2
         \right)^{1/2}+ \left(\frac{1}{T} \sum_{t=1}^T \norm{\H_1\f_t-\hat{\f}_t^{(1)}}^2 \cdot \frac{1}{T} \sum_{t=1}^T
         \norm{\hat{\f}_t^{(1)}}^2 \right)^{1/2} \\
       &=O_P\left(\frateone \right),
  \end{align*}
  where the last equality follows from Lemma \ref{lem:A5} and that
  $\norm{\H_1\f_t}\leq \norm{\H_1} \norm{\f_t}=O_P\left(1\right)$ and
  $\norm{\hat{\f}_t^{(1)}}=O_P\left(1\right)$. Similarly, $\fnorm{\H_2\H_2'-\I_K}=O_P\left(\fratetwo\right)$.

  (iii) Since $\fnorm{\H_1\H_1'-\I_K}=O_P\left(\frateone\right)$ and $\norm{\H_1}=O_P\left(1\right)$, we have
  $\fnorm{\H_1\H_1'\H_1-\H_1}=O_P\left(\frateone\right)$. Since $\H_1^{-1}=\H_1^{-1}(\I_K-\H_1\H_1'+\H_1\H_1')$, it
  follows Lemma \ref{lem:A3}(i) that $\norm{\H_1^{-1}}\leq \norm{\H_1^{-1}}O_P\left(\frateone\right)+\norm{\H_1'}$. Hence,
  $\norm{\H_1^{-1}}=O_P\left(1\right)$. Left multiplying $\H_1\H_1'\H_1-\H_1$ by $\H_1^{-1}$ gives
  $\fnorm{\H_1'\H_1-\I_K}=O_P\left(\frateone\right)$. Similarly, $\fnorm{\H_2'\H_2-\I_K}=O_P\left(\fratetwo\right)$.
\end{proof}

\begin{lem}
  \label{lem:A4}
  Let $\C_1=\hat{\B}_1-\B_S\H_1'$ and $\C_2=\hat{\B}_2-\B_S\H_2'$, where
  $\hat{\B}_1$, $\hat{\B}_2$, and $\B_S$ are defined in Section \ref{sec:4}.\\
  (i) $\fnorm{\C_1}^2= O_P\left(sw_1^2\right)$, $\fnorm{\C_2}^2= O_P\left(sw_2^2\right)$;
  $\snorm{\C_1\C_1'}^2=O_P\left(sw_1^4\right)$,
  $\snorm{\C_2\C_2'}^2=O_P\left(s w_2^4\right)$.\\
  (ii) $\snorm{\hat{\Sigmab}_{u,S}^{(1)}-\Sigmab_{u,S}}^2=O_P\left(m_s^2w_1^2\right)$;
  $\snorm{\hat{\Sigmab}_{u,S}^{(2)}-\Sigmab_{u,S}}^2=O_P\left(m_s^2w_2^2\right)$.\\
  (iii) $\snorm{\B_S\H_1'\C_1'}^2=O_P\left(w_1^2\right)$;
  $\snorm{\B_S\H_2'\C_2'}^2=O_P\left(w_2^2\right)$.\\
  (iv) $\snorm{\B_S(\H_1'\H_1-\I_K)\B_S'}^2=O_P\left(\frac{1}{sT}+\frac{1}{s^2} \right)$;
  $\snorm{\B_S(\H_2'\H_2-\I_K)\B_S'}^2=O_P\left(\frac{1}{sT}+\frac{1}{sp} \right)$.
\end{lem}
\begin{proof}[\textbf{Proof of Lemma \ref{lem:A4}.}]
  (i) We have
  $\fnorm{\C_1}^2\leq s(\max_{i\leq s} \norm{\hat{\b}^{(1)}_i-\H\b_i})^2=O_P\left( sw_1^2 \right)$. By the
  general result that for any matrix $\A$,
  $\snorm{\A}^2=s^{-1}\fnorm{\Sigmab_S^{-1/2}\A\Sigmab_S^{-1/2}}^2=O_P\left(s^{-1} \fnorm{\A}^2 \right)$, we
  have $\snorm{\C_1'\C_1}^2=O_P\left(s^{-1} \fnorm{\C_1}^4 \right)=O_P\left( s w_1^4\right)$.  Similarly,
  $\fnorm{\C_2}^2=O_P\left(s w_2^2 \right)$ and $\snorm{\C_2\C_2'}^2=O_P\left(sw_2^4 \right)$.

  (ii) By Lemma \ref{lem:3},
  \begin{equation*}
    \snorm{\hat{\Sigmab}_{u,S}^{(1)}-\Sigmab_{u,S}}^2=O_P\left(s^{-1} \fnorm{\hat{\Sigmab}_{u,S}^{(1)}-\Sigmab_{u,S}}^2 \right)=
    O_P\left(\norm{\hat{\Sigmab}_{u,S}^{(1)}-\Sigmab_{u,S}}^2
    \right)=O_P\left(m_s^2w_1^2\right).
  \end{equation*}
  Similar results can be shown for
  $\snorm{\hat{\Sigmab}_{u,S}^{(2)}-\Sigmab_{u,S}}$.

  (iii) By adapt the proof of Theorem 2 in \cite{FFL08}, we have that
  $\norm{\B_S'\Sigmab_S^{-1}\B_S}=O(1)$. Hence,
  \begin{align*}
    \snorm{\B_S\H_1'\C_1'}^2
    &= s^{-1}
      \tr(\H_1'\C_1'\Sigmab_S^{-1}\C_1\H_1\B_S'\Sigmab_S^{-1}\B_S)\\
    &\leq s^{-1} \norm{\H_1}^2 \norm{\B_S'\Sigmab_S^{-1}\B_S}
      \norm{\Sigmab_S^{-1}} \fnorm{\C_1}^2\\
    &=O_P\left(s^{-1} \fnorm{\C_1}^2 \right)=O_P\left(w_1^2\right).
  \end{align*}
  Similarly, $\snorm{\B_S\H_2'\C_2'}=O_P\left(w_2^2 \right)$.

  (iv) We have
  \begin{align*}
    \snorm{\B_S(\H_1'\H_1-\I_K)\B_S'}^2
    &= s^{-1}
      \tr((\H_1'\H_1-\I_K)\B_S'\Sigmab_S^{-1}\B_S(\H_1'\H_1-\I_K)\B_S'\Sigmab_S^{-1}\B_S) \\
    &\leq s^{-1} \fnorm{\H_1'\H_1-\I_K}^2
      \norm{\B_S'\Sigmab_S^{-1}\B_S}^2 =
      O_P\left(\frac{1}{sT}+\frac{1}{s^2}\right).
  \end{align*}
  Similarly, $\snorm{\B_S(\H_2'\H_2-\I_K)\B_S'}^2=O_P\left(\frac{1}{sT}+\frac{1}{sp} \right)$.
\end{proof}

\begin{lem}
  \label{lem:A5}
  Under conditions of Lemma \ref{lem:1},
  \begin{align*}
    \frac{1}{T} \sum_{t=1}^T \norm{\hat{\f}_t^{(1)}-\H_1\f_t}^2&=O_P\left(1/s+1/T \right)\\
    \frac{1}{T} \sum_{t=1}^T \norm{\hat{\f}_t^{(2)}-\H_2\f_t}^2&=O_P\left(1/p+1/T \right)
  \end{align*}
\end{lem}
\begin{proof}[\textbf{Proof of Lemma \ref{lem:A5}.}]
  Without loss of generality, we only prove the result for general $p$. Again, we
  write $\hat{\f}_t^{(2)}$ as $\hat{\f}_t$, $\H_2$ as $\H$ and $\hat{\V}_2$ as
  $\hat{\V}$ for notational simplicity. By (\ref{eq:A2}),
  \begin{equation*}
    \frac{1}{T} \sum_{t=1}^T \norm{\hat{\f}_t-\H\f_t}^2 \leq c \norm{\hat{\V}^{-1}}^2 \sum_{j=1}^7 \frac{1}{T}
    \sum_{t=1}^T G_{jt}^2,
  \end{equation*}
  where $c$ is a positive constant and $G_{jt}$ is the $j$th summand on the right hand side of (\ref{eq:A2}). By Lemma A.6
  (iv) of \cite{BW13}, $\frac{1}{T} \sum_{i=1}^T G_{1t}^2=o_P\left(1/p+1/T \right)$. By Lemma A.10 (i) and (iii) of
  \cite{BW13}, $\frac{1}{T} \sum_{t=1}^T G_{2t}^2=O_P\left(1/T \right)$ and
  $\frac{1}{T} \sum_{t=1}^T G_{3t}^2=O_P\left(1/T \right)$.  By Lemma A.6 (iii), (v) and (vi) of \cite{BW13},
  $\frac{1}{T} \sum_{t=1}^T G_{4t}^2=o_P(1/p)$, $\frac{1}{T} \sum_{t=1}^T G_{6t}^2=o_P(1/p)$ and
  $\frac{1}{T} \sum_{t=1}^T G_{7t}^2=o_P(1/p)$. Finally, by Lemma A.11 (ii) of \cite{BW13},
  $\frac{1}{T} \sum_{t=1}^T G_{5t}^2=O_P(1/p)$. Therefore, the dominating terms are $G_{2t}$, $G_{3t}$ and $G_{5t}$, which
  together give the rate of $O_P\left(1/p+1/T \right)$.
\end{proof}

\begin{lem}
  \label{lem:A6}
  With probability tending to 1,\\
  (i) $\lambda_{\min}(\I_K+(\B_S\H_1')'\Sigmab_{u,S}^{-1}\B_S\H_1')\geq cs$,
  $\lambda_{\min}(\I_K+(\B_S\H_2')'\Sigmab_{u,S}^{-1}\B_S\H_2')\geq
  cs$, $\lambda_{\min}(\I_K+\B_S'\Sigmab_{u,S}^{-1}\B_S)\geq cs$; \\
  (ii) $\lambda_{\min} (\I_K+\hat{\B}_1'(\hsigus{(1)})^{-1}\hat{\B}_1)\geq cs$,
  $\lambda_{\min} (\I_K+\hat{\B}_2'(\hsigus{(2)})^{-1}\hat{\B}_2)\geq cs$,
  $\lambda_{\min} (\I_K+\hat{\B}_o'(\hsigus{o})^{-1}\hat{\B}_o)\allowbreak\geq cs$;\\
  (iii) $\lambda_{\min} ((\H_1'\H_1)^{-1}+\B_S'\Sigmab_{u,S}^{-1}\B_S) \geq cs$,
  $\lambda_{\min} ((\H_2'\H_2)^{-1}+\B_S'\Sigmab_{u,S}^{-1}\B_S) \geq cs$.
\end{lem}
\begin{proof}[\textbf{Proof of Lemma \ref{lem:A6}.}]
  By Lemma \ref{lem:A3}, with probability tending to one, $\lambda_{\min}(\H_1\H_1')$ is bounded away from 0. Therefore,
  \begin{align*}
    &\hspace{3ex}\lambda_{\min}(\I_K+(\B_S\H_1')' \Sigmab_{u,S}^{-1} \B_S\H_1')\geq
      \lambda_{\min}(\H_1\B'_S\Sigmab_{u,S}^{-1}\B_S\H_1') \\
    &\geq
      \lambda_{\min}(\Sigmab_{u,S}^{-1})\lambda_{\min}(\B'_S\B_S)\lambda_{\min}(\H_1\H_1') \geq cs.
  \end{align*}
  Similar results hold for the other two statements. The results in (ii) follow from (i) and (\ref{eq:A6}).  The statement
  (iii) follows from a similar argument as $\H_1\H'_1$ and $\H_2\H'_2$ are positive semi-definite.
\end{proof}

\begin{lem}
  \label{lem:A7} [Theorem A.1 of \cite{FLM13}] 
  Let $\hat{u}_{it}$ be defined as in step ii. of Method 1 in Section \ref{sec:4}. Under conditions (iv),
  (v), if there is a sequence $a_T=o(1)$ so that
  $\max_{i\leq p} \frac{1}{T} \sum_{t=1}^T |u_{it}-\hat{u}_{it}|^2=O_P\left(a_T^2 \right)$ and
  $\max_{i\leq p, t\leq T} |u_{it}-\hat{u}_{it}|=o_P\left(1 \right)$, then the adaptive thresholding estimator
  $\hat{\Sigmab}_u$ with $\omega(p)=\sqrt{(\log p)/T}+a_T$ satisfies that
  $\norm{\hat{\Sigmab}_u-\Sigmab_u}=O_P\left(m_p[\omega(p)]^{1-q} \right)$. If further
  $m_p[\omega(p)]^{1-q}=o(1)$, then $\hat{\Sigmab}_u$ is invertible with probability approaching one, and
  $\norm{\hat{\Sigmab}^{-1}_u-\Sigmab^{-1}_u}=O_P\left(m_p[\omega(p)]^{1-q} \right)$.
\end{lem}

\end{document}